\documentclass[a4, 12pt, dvipdfmx, leqno]{article}
\usepackage{amsmath}
\usepackage{amssymb}
\usepackage{graphicx}
\usepackage{overpic}
\usepackage{here}
\newtheorem{df}{Definition}[section]

\newtheorem{Th}[df]{Theorem}
\newtheorem{Rm}[df]{Remark}
\newtheorem{lm}[df]{Lemma}

\newcommand{\artanh}{\mathrm{artanh \hspace*{0.6mm}}}

\begin{document}

\title{\bf On the Existence of 
Leapfrogging Pair of Circular Vortex Filaments
}
\author{Masashi A{\sc iki}}
\date{}
\maketitle
\vspace*{-0.5cm}

\begin{abstract}
We propose and analyze a system of nonlinear partial differential equations
describing the motion of a pair of vortex filaments. Furthermore, 
for a pair of coaxial circular vortex filaments,
we derive a condition for leapfrogging to occur and 
prove that the condition is
 necessary and sufficient for the 
occurrence of leapfrogging.
\end{abstract}

\section{Introduction and Problem Setting}

In this paper, we are interested 
in the interaction of two vortex rings sharing the same axis of 
symmetry. The study of the interaction of two vortex rings 
dates back to the classical paper by Helmholtz \cite{29},
where he observed that a 
pair of vortex rings may exhibit what is now known as 
``leapfrogging''. Leapfrogging is a motion pattern where 
two vortex rings sharing a common axis of symmetry pass through
each other repeatedly due to the induced flow of the rings acting 
on each other. 
Dyson \cite{27,28} considered the motion of 
vortex rings and
proposed a system of equations describing the motion of 
a pair of coaxial vortex rings. Based on this model system, 
Dyson also observed that leapfrogging may occur.
Hicks \cite{30} also considered the interaction of a pair of 
vortex rings and derived a model similar to the one obtained by 
Dyson. He also made numerical observations that depending on the vorticity strengths and 
initial configuration of the rings, the pair may show leapfrogging, 
or the two rings may separate indefinitely. 
Although these observations were known for a long time, 
the first experiment which successfully provided 
photographic proofs of
the leapfrogging phenomenon in a laboratory setting was 
conducted by Yamada and Matsui \cite{32}. 
They used vortex rings made of air and used smoke to visualize the 
rings and created a leapfrogging pair of rings, and hence, 
leapfrogging vortex rings were observed in the real world. 
In recent years, detailed models 
for the motion of vortex rings have been obtained,
for example by Saffman \cite{41} and Fukumoto \cite{42},
in which various distributions of vorticity in the ring core
and the change in shape of the core can be incorporated.
Based on these models, Borisov, Kilin, and Mamaev \cite{39}
gave a complete description of the possible motion patterns 
of two interacting vortex rings with a common axis of symmetry.
Shariff and Leonard \cite{40} give a review of the 
history of the research of vortex rings
in which many other phenomena related to
the motion of single and multiple 
vortex rings, including the 
leapfrogging phenomenon, are addressed.

The study of the motion of a vortex ring, and the interaction of vortex rings are
not only interesting from a theoretical stand point, 
but also plays an important role in engineering.
One such example being the generation of sound from a round jet as studied in 
Crighton \cite{33}, Hussain and Zaman \cite{34}, and Zaman \cite{35}.

The aim of the present paper is to propose and 
analyze a system of nonlinear partial differential equations
for the motion of two interacting vortex filaments.
A vortex filament is 
a space curve on which the vorticity of the fluid is concentrated.
A vortex filament can be seen as an idealization of a thin vortex 
structure for which the 
evolution can be approximated by the evolution of its center line.
A model equation,
known as the Localized Induction Equation, describing the motion of a 
single vortex filament was first proposed by 
Da Rios \cite{20} and later independently derived by 
Murakami et al. \cite{22} and Arms and Hama \cite{21}.
The main idea used in these papers is the so-called 
localized induction approximation, and we follow this concept in this 
paper to derive a system of equations describing the motion of 
two interacting vortex filaments. We further consider  
the case when the filaments are circular with a common
axis of symmetry and obtain a 
necessary and sufficient condition for leapfrogging to occur. 

The motivation for this work is the following.
Many of the analysis made for leapfrogging vortex rings 
are conducted for a system of ordinary differential equations. 
This, of course, is natural since
for a coaxial vortex ring pair, the circular shape is 
expected to be preserved from the symmetry of the induced flow. On the other hand, when one considers
stability of such motion, it would be useful if the motion is 
described in the framework of partial differential equations because
it becomes possible to consider non-symmetric perturbations.
Consequently, this complicates the situation and thus, we 
consider circular vortex filaments instead of 
vortex rings with finite core thickness. 
As far as the author knows, the model
proposed by Klein, Majda, and Damodaran \cite{38} is the only 
model describing the motion of multiple vortex filaments in the framework of 
partial differential equations.
They consider filaments that are nearly straight and parallel to each other 
and derive a model system describing their motion, and as such, 
the motion described by this model is two-dimensional.
Hence, when considering leapfrogging circular filaments, 
the model in \cite{38} is not suitable, and we 
derive a different system in this paper.

The rest of the paper is organized as follows.
In Section 2, we derive the model system  
via the localized induction approximation. 
We also give some exact solutions of the obtained system to 
show that the model is capable of describing well known
motions of straight vortex filaments which are parallel to each other.
In Section 3, we consider the case when the two filaments 
are circular with a common axis of symmetry and 
the vorticity strengths have the same sign. We show that 
the problem can be reduced to a two-dimensional Hamiltonian system
and give a necessary and sufficient condition for 
leapfrogging to occur. The precise statement will be 
given in the beginning of Section 3. 
In Section 4, the leapfrogging phenomenon will be considered for 
a pair of circular filaments which have vorticity strengths of opposite signs.
Similar to Section 3, we give a necessary and sufficient condition for leapfrogging to
occur. Finally in Section 5, 
we compare our results with the results obtained in 
Borisov, Kilin, and Mamaev \cite{39}, and also give 
concluding remarks.


\section{Interaction of Two Vortex Filaments}
\setcounter{equation}{0}
We consider the interaction of two vortex filaments and 
derive a system of nonlinear partial differential equations
which describe their motion.
The obtained model admits solutions which correspond to 
well known motions of point vortices when the two filaments 
are straight parallel lines, 
and also gives a clear view of the dynamics when the filaments are arranged as 
coaxial circles, 
and hence the author hopes 
that the model could be of some significance.
%
%
%
%
%
%
\subsection{Derivation of the Model System}
Following the work of Arms and Hama \cite{21}, we 
apply the localized induction approximation to the 
Biot--Savart law to obtain a system of partial differential equations
approximating the motion of two interacting vortex filaments.
The velocity
\( \mbox{\mathversion{bold}$v$}(\mbox{\mathversion{bold}$x$})\)
at some point \( \mbox{\mathversion{bold}$x$}\in \mathbf{R}^{3}\)
of an infinite body of
incompressible and inviscid fluid induced by a pair of vortex filaments 
whose positions are parametrized by \( \xi\in J \) at time 
\( t\geq 0 \) as \( \mbox{\mathversion{bold}$X$}(\xi ,t)\) and 
\( \mbox{\mathversion{bold}$Y$}(\xi ,t)\) is 
given by
\begin{align}
\mbox{\mathversion{bold}$v$}
(\mbox{\mathversion{bold}$x$})&=
\frac{\Gamma_{1}}{4\pi}\int_{J}
\frac{\mbox{\mathversion{bold}$X$}_{\xi}(r,t)\times
(\mbox{\mathversion{bold}$x$}-\mbox{\mathversion{bold}$X$}(r,t))}
{|\mbox{\mathversion{bold}$x$}-\mbox{\mathversion{bold}$X$}(r,t)|^{3}}
\ {\rm d}r
+
\frac{\Gamma_{2}}{4\pi}
\int_{J}
\frac{\mbox{\mathversion{bold}$Y$}_{\xi}(r,t)\times
(\mbox{\mathversion{bold}$x$}-\mbox{\mathversion{bold}$Y$}(r,t))}
{|\mbox{\mathversion{bold}$x$}-\mbox{\mathversion{bold}$Y$}(r,t)|^{3}}
\ {\rm d}r 
\label{biot}
\end{align}
where \( \times \) is the exterior product in the three-dimensional 
Euclidean space, \( \Gamma_{1}\) is the vorticity strength of 
the filament \( \mbox{\mathversion{bold}$X$}\),
\( \Gamma_{2}\) is the vorticity strength of the filament 
\( \mbox{\mathversion{bold}$Y$} \), \( J=\mathbf{R} \ \mbox{or} \ \mathbf{R}/2\pi \mathbf{Z}\), and subscripts denote
the partial differentiation with the respective variables.
The above equation is the Biot--Savart law when the 
vorticity is concentrated on two vortex filaments. 
The case \( J=\mathbf{R}\) corresponds to when \( \mbox{\mathversion{bold}$X$}\) and 
\( \mbox{\mathversion{bold}$Y$}\) are infinitely long filaments, and 
the case \( J=\mathbf{R}/2\pi \mathbf{Z}\) corresponds to when 
\( \mbox{\mathversion{bold}$X$}\) and \( \mbox{\mathversion{bold}$Y$}\) are closed filaments. 
To determine the velocity of a point on one of the filaments
(say \( \mbox{\mathversion{bold}$X$}(\xi ,t)\)), one would like to
substitute \( \mbox{\mathversion{bold}$x$}=
\mbox{\mathversion{bold}$X$}(\xi ,t)\) in equation 
(\ref{biot}), but this would result in the divergence of the 
first integral on the right-hand side. Hence we apply the localized 
induction approximation to approximate the the velocity 
at \( \mbox{\mathversion{bold}$X$}(\xi ,t)\) by the following 
equation.
\begin{align*}
\mbox{\mathversion{bold}$v$}
(\mbox{\mathversion{bold}$X$}(\xi ,t))&=
\frac{\Gamma_{1}}{4\pi}
\int_{\varepsilon \leq |\xi-r| \leq L}
\frac{\mbox{\mathversion{bold}$X$}_{\xi}(r,t)\times
(\mbox{\mathversion{bold}$X$}(\xi ,t)
-\mbox{\mathversion{bold}$X$}(r,t))}
{|\mbox{\mathversion{bold}$X$}(\xi ,t)
-\mbox{\mathversion{bold}$X$}(r,t)|^{3}}
\ {\rm d}r \\[3mm]
& \qquad +
\frac{\Gamma_{2}}{4\pi}
\int_{|\xi-r|\leq \delta  }
\frac{\mbox{\mathversion{bold}$Y$}_{\xi}(r,t)\times
(\mbox{\mathversion{bold}$X$}(\xi ,t)
-\mbox{\mathversion{bold}$Y$}(r,t))}
{|\mbox{\mathversion{bold}$X$}(\xi ,t)
-\mbox{\mathversion{bold}$Y$}(r,t)|^{3}}
\ {\rm d}r \\[3mm]
&=: I_{1}+I_{2}.
\end{align*}
Here, \( \varepsilon >0 \) and \( \delta >0\) are
small parameters, and \( L>0\) is a cut-off parameter. \( I_{1}\) is the effect of self-induction, 
and \( I_{2}\) is the effect of interaction.
The approximation applied in \( I_{1}\) is the well known 
localized induction approximation.
To obtain \( I_{2}\), we have further assumed that 
the filaments \( \mbox{\mathversion{bold}$X$}\) and \( \mbox{\mathversion{bold}$Y$}\) are 
positioned in a way that 
\( \mbox{\mathversion{bold}$Y$}(\xi ,t)\) is the 
closest point to 
\( \mbox{\mathversion{bold}$X$}(\xi ,t)\)
and the 
contributions from points far away from 
\( \mbox{\mathversion{bold}$Y$}(\xi ,t)\) can be ignored.
This kind of geometric assumption is true for the 
situations that we treat in this paper, 
but does not hold, for example, when the filaments are 
knotted together.
By the calculations 
in Arms and Hama \cite{21}, it is known that 
\( I_{1}\) can be expanded in terms of small \( \varepsilon \)
as follows.
\begin{align*}
I_{1}=-\frac{\Gamma_{1}}{4\pi}\log(\frac{L}{\varepsilon})
\frac{\mbox{\mathversion{bold}$X$}_{\xi}\times
\mbox{\mathversion{bold}$X$}_{\xi \xi }}
{|\mbox{\mathversion{bold}$X$}_{\xi }|^{3}}
+
O(1).
\end{align*}
The above is obtained by substituting the 
Taylor expansion of 
\( \mbox{\mathversion{bold}$X$}(r,t) \) and 
\( \mbox{\mathversion{bold}$X$}_{\xi}(r,t) \)
with respect to \( r\) around \( \xi \) into the 
integrand. We further substitute
\begin{align*}
\mbox{\mathversion{bold}$Y$}(r,t)&=
\mbox{\mathversion{bold}$Y$}(\xi ,t) + 
\mbox{\mathversion{bold}$Y$}_{\xi}(\xi ,t)(r-\xi)
+
O((r-\xi)^{2}), \\[3mm]
\mbox{\mathversion{bold}$Y$}_{\xi}(r,t)&=
\mbox{\mathversion{bold}$Y$}_{\xi}(\xi, t) +
\mbox{\mathversion{bold}$Y$}_{\xi \xi}(\xi ,t)(r-\xi) +
O((r-\xi)^{2}),
\end{align*}
into \( I_{2}\) to obtain 
\begin{align*}
I_{2}=
\frac{\delta \Gamma_{2}}{2\pi}
\frac{\mbox{\mathversion{bold}$Y$}_{\xi}\times
(\mbox{\mathversion{bold}$X$}-\mbox{\mathversion{bold}$Y$})}
{|\mbox{\mathversion{bold}$X$}-\mbox{\mathversion{bold}$Y$}|^{3}}
+
O(\delta^{2}).
\end{align*}
Hence, after fixing \( L\) and taking sufficiently small \( \varepsilon\) 
and \( \delta \), the leading order terms of \( I_{1}\) and 
\( I_{2}\) yield
\begin{align*}
\mbox{\mathversion{bold}$X$}_{t}
=
-\frac{\Gamma_{1}}{4\pi}\log(\frac{L}{\varepsilon})
\frac{\mbox{\mathversion{bold}$X$}_{\xi}\times
\mbox{\mathversion{bold}$X$}_{\xi \xi }}
{|\mbox{\mathversion{bold}$X$}_{\xi }|^{3}}
+
\frac{\delta \Gamma_{2}}{2\pi}
\frac{\mbox{\mathversion{bold}$Y$}_{\xi}\times
(\mbox{\mathversion{bold}$X$}-\mbox{\mathversion{bold}$Y$})}
{|\mbox{\mathversion{bold}$X$}-\mbox{\mathversion{bold}$Y$}|^{3}},
\end{align*}
where we also used the fact that
\( \mbox{\mathversion{bold}$v$}
(\mbox{\mathversion{bold}$X$}(\xi ,t))=
\mbox{\mathversion{bold}$X$}_{t}(\xi ,t) \) by the definition of 
velocity. 
By rescaling time by a factor of 
\( -\log (\frac{L}{\varepsilon})/4\pi \), we obtain
\begin{align*}
\mbox{\mathversion{bold}$X$}_{t}=
\Gamma_{1}\frac{
\mbox{\mathversion{bold}$X$}_{\xi}\times 
\mbox{\mathversion{bold}$X$}_{\xi \xi}}
{
|\mbox{\mathversion{bold}$X$}_{\xi}|^{3}
}
-
\alpha \Gamma_{2}
\frac{\mbox{\mathversion{bold}$Y$}_{\xi}\times
(\mbox{\mathversion{bold}$X$}-\mbox{\mathversion{bold}$Y$})}
{|\mbox{\mathversion{bold}$X$}-\mbox{\mathversion{bold}$Y$}|^{3}},
\end{align*}
where \( \alpha = 2\delta /\log (\frac{L}{\varepsilon})>0 \). 
The calculations for the 
velocity at points on \( \mbox{\mathversion{bold}$Y$}\) are the same
and hence we arrive at the following system.
\begin{align}
\left\{
\begin{array}{l}
\displaystyle
\mbox{\mathversion{bold}$X$}_{t}=
\Gamma_{1}\frac{
\mbox{\mathversion{bold}$X$}_{\xi}\times 
\mbox{\mathversion{bold}$X$}_{\xi \xi}}
{
|\mbox{\mathversion{bold}$X$}_{\xi}|^{3}
}
-
\alpha \Gamma_{2}
\frac{\mbox{\mathversion{bold}$Y$}_{\xi}\times
(\mbox{\mathversion{bold}$X$}-\mbox{\mathversion{bold}$Y$})}
{|\mbox{\mathversion{bold}$X$}-\mbox{\mathversion{bold}$Y$}|^{3}},
\\[5mm]
\displaystyle
\mbox{\mathversion{bold}$Y$}_{t}=
\Gamma_{2}\frac{
\mbox{\mathversion{bold}$Y$}_{\xi}\times 
\mbox{\mathversion{bold}$Y$}_{\xi \xi}}
{
|\mbox{\mathversion{bold}$Y$}_{\xi}|^{3}
}
-
\alpha \Gamma_{1}
\frac{\mbox{\mathversion{bold}$X$}_{\xi}\times
(\mbox{\mathversion{bold}$Y$}-\mbox{\mathversion{bold}$X$})}
{|\mbox{\mathversion{bold}$X$}-\mbox{\mathversion{bold}$Y$}|^{3}}.
\end{array}\right.
\label{model}
\end{align}
All the analysis that follows will be based on the above system (\ref{model}).
%
%
%
%
\subsection{Dynamics of Two Parallel Lines}
As a preliminary analysis, we show that for a pair of infinitely long, straight, and parallel vortex filaments, 
the dynamics of the filaments according to equation (\ref{model}) are the same as that of two point vortices moving in a plane.
Suppose the two filaments are initially parametrized as
\begin{align*}
\mbox{\mathversion{bold}$X$}_{0}(\xi )=
{}^{t}(l,0,\xi), \qquad 
\mbox{\mathversion{bold}$Y$}_{0}(\xi )=
{}^{t}(-l,0,\xi ),
\end{align*}
where \( l>0\) is arbitrary. In this situation, it is expected that the motions of the filaments become two-dimensional and 
resemble that of two point vortices. Indeed, if we make the ansatz
\begin{align*}
\mbox{\mathversion{bold}$X$}(\xi ,t)=
{}^{t}(x_{1}(t),x_{2}(t),\xi ), \qquad 
\mbox{\mathversion{bold}$Y$}(\xi ,t)=
{}^{t}(y_{1}(t),y_{2}(t),\xi),
\end{align*}
and substitute it into (\ref{model}), we obtain
\begin{align*}
\left\{
\begin{array}{l}
\displaystyle
\dot{x_{1}}= \frac{\alpha \Gamma_{2}(x_{2}-y_{2})}
{\big( (x_{1}-x_{1})^{2}+(x_{2}-y_{2})^{2}\big)^{3/2}}, \\[5mm]
\displaystyle 
\dot{x_{2}}= -\frac{\alpha \Gamma_{2}(x_{1}-y_{1})}
{\big( (x_{1}-y_{1})^{2}+(x_{2}-y_{2})^{2}\big)^{3/2}}, \\[5mm]
\displaystyle 
\dot{y_{1}}= \frac{\alpha \Gamma_{1}(y_{2}-x_{2})}
{\big( (x_{1}-x_{1})^{2}+(x_{2}-y_{2})^{2}\big)^{3/2}}, \\[5mm]
\displaystyle
\dot{y_{2}}=  -\frac{\alpha \Gamma_{1}(y_{1}-x_{1})}
{\big( (x_{1}-x_{1})^{2}+(x_{2}-y_{2})^{2}\big)^{3/2}}, \\[5mm]
\end{array}\right.
\end{align*}
where a dot over a variable denotes the derivative with respect to time.
Further setting \( z_{1}=x_{1}+ix_{2}\) and \( z_{2}=y_{1}+iy_{2}\), where 
\( i\) is the imaginary unit, we have
\begin{align*}
\left\{
\begin{array}{l}
\displaystyle
\dot{z_{1}}=-i\alpha \Gamma_{2}\frac{z_{1}-z_{2}}{|z_{1}-z_{2}|^{3/2}}, \\[5mm]
\displaystyle
\dot{z_{2}}=-i\alpha \Gamma_{1}\frac{z_{2}-z_{1}}{|z_{1}-z_{2}|^{3/2}}.
\end{array}\right.
\end{align*}
We see from direct calculation that when \( \Gamma_{1}+\Gamma_{2}\neq 0\),
\begin{align*}
C = \frac{\Gamma_{1}z_{1}+\Gamma_{2}z_{2}}{\Gamma_{1}+\Gamma_{2}}, \qquad 
D = |z_{1}-z_{2}|,
\end{align*}
are conserved quantities. \( C\) is known as the center of vorticity.
Utilizing these quantities, the equations can be decoupled to obtain
\begin{align*}
\begin{pmatrix}
\dot{z_{1}}\\[3mm]
\dot{z_{2}}
\end{pmatrix}
=
-\frac{i\alpha(\Gamma_{1}+\Gamma_{2})}{D^{3/2}}
\begin{pmatrix}
1 &  \ 0 \\[3mm]
0 &  \ 1 
\end{pmatrix}
\begin{pmatrix}
z_{1}-C\\[3mm]
z_{2}-C
\end{pmatrix}.
\end{align*}
The above equation can be solved explicitly and we have
\begin{align*}
z_{j}(t)=(z_{j}(0)-C)e^{i\omega t}+C
\end{align*}
for \( j=1,2\), where 
\( \omega = -\alpha (\Gamma_{1}+\Gamma_{2})/D^{3/2}\).
This shows that the two filaments rotate in a two-dimensional circular pattern and the 
center and radius of rotation is determined by the center of vorticity.
When \( \Gamma_{1}+\Gamma_{2}=0\), we see that 
\( z_{1}-z_{2} \) is conserved and hence we have
\begin{align*}
\dot{z_{j}}=-\frac{i\alpha \Gamma_{2}}{D^{3/2}}w_{0} = \mbox{const.},
\end{align*}
for \( j=1,2\) with \( w_{0}=z_{1}(0)-z_{2}(0)\). This shows that the two filaments 
travel in a straight line at a constant speed while keeping their parallel configuration.
These dynamics of the filaments directly correspond to the motion of two point vortices
moving in a plane, which is well known in the literature such as Newton \cite{36}.
Hence, we see that system (\ref{model}) is capable of describing 
the motion of two parallel lines in the expected manner.
%
%
%
%
%
\section{Leapfrogging for a Pair of Filaments with Vorticity Strengths of 
the Same Sign}
\setcounter{equation}{0}
We consider the case when the two filaments are arranged as coaxial circles and \( \Gamma_{1},\Gamma_{2}>0\). 
Rescaling the time variable by a factor of \( \Gamma_{2}\)
in (\ref{model}) yields
\begin{align}
\left\{
\begin{array}{l}
\displaystyle
\mbox{\mathversion{bold}$X$}_{t}=
\beta \frac{
\mbox{\mathversion{bold}$X$}_{\xi}\times 
\mbox{\mathversion{bold}$X$}_{\xi \xi}}
{
|\mbox{\mathversion{bold}$X$}_{\xi}|^{3}
}
-
\alpha 
\frac{\mbox{\mathversion{bold}$Y$}_{\xi}\times
(\mbox{\mathversion{bold}$X$}-\mbox{\mathversion{bold}$Y$})}
{|\mbox{\mathversion{bold}$X$}-\mbox{\mathversion{bold}$Y$}|^{3}},
\\[5mm]
\displaystyle
\mbox{\mathversion{bold}$Y$}_{t}=
\frac{
\mbox{\mathversion{bold}$Y$}_{\xi}\times 
\mbox{\mathversion{bold}$Y$}_{\xi \xi}}
{
|\mbox{\mathversion{bold}$Y$}_{\xi}|^{3}
}
-
\alpha \beta 
\frac{\mbox{\mathversion{bold}$X$}_{\xi}\times
(\mbox{\mathversion{bold}$Y$}-\mbox{\mathversion{bold}$X$})}
{|\mbox{\mathversion{bold}$X$}-\mbox{\mathversion{bold}$Y$}|^{3}},
\end{array}\right.
\label{model2}
\end{align}
where \( \beta = \Gamma_{1}/\Gamma_{2}\).
We assume without loss of generality that \( \beta \geq 1\),
since the case \( \beta <1 \) is reduced to the case
\( \beta>1\) by renaming the filaments.

Suppose that for some \( R_{1,0},R_{2,0}>0\) and \( z_{1,0},z_{2,0}\in \mathbf{R}\), 
the initial filaments 
\( \mbox{\mathversion{bold}$X$}_{0}\) and \( \mbox{\mathversion{bold}$Y$}_{0}\)
are parametrized by \( \xi \in [0,2\pi )\) as follows.
\begin{align*}
\mbox{\mathversion{bold}$X$}_{0}(\xi )={}^{t}
(R_{1,0}\cos(\xi), R_{1,0}\sin(\xi), z_{1,0}), \quad 
\mbox{\mathversion{bold}$Y$}_{0}(\xi)={}^{t}
(R_{2,0}\cos(\xi), R_{2,0}\sin (\xi), z_{2,0}),
\end{align*}
where we assume that \( (R_{1,0}-R_{2,0})^{2}+(z_{1,0}-z_{2,0})^{2}>0\), 
which means that the two circles are not overlapping. 
Now, we make the ansatz
\begin{align*}
\mbox{\mathversion{bold}$X$}(\xi ,t)={}^{t}
(R_{1}(t)\cos (\xi), R_{1}(t)\sin(\xi), z_{1}(t)), \quad 
\mbox{\mathversion{bold}$Y$}(\xi ,t)={}^{t}
(R_{2}(t)\cos (\xi), R_{2}(t)\sin (\xi), z_{2}(t)),
\end{align*}
and substitute it into (\ref{model2}).
From the equation for \( \mbox{\mathversion{bold}$X$}\) we have
\begin{align*}
\dot{R_{1}}\cos (\xi)&=-\frac{\alpha R_{2}(z_{1}-z_{2})\cos (\xi)}{\big( (R_{1}-R_{2})^{2}+(z_{1}-z_{2})^{2}\big)^{3/2}},\\[5mm]
\dot{R_{1}}\sin (\xi)&=-\frac{\alpha R_{2}(z_{1}-z_{2})\sin (\xi)}{\big( (R_{1}-R_{2})^{2}+(z_{1}-z_{2})^{2}\big)^{3/2}},\\[5mm]
\dot{z_{1}}&=\frac{\beta}{R_{1}}+\frac{\alpha R_{2}(R_{1}-R_{2})}{\big( (R_{1}-R_{2})^{2}+(z_{1}-z_{2})^{2}\big)^{3/2}}.
\end{align*}
The dependence of the system on \( \xi \) is eliminated by multiplying the first two equations by
\( \cos (\xi) \) and \( \sin (\xi)\), respectively, and adding. The equations for
\( \mbox{\mathversion{bold}$Y$}\) are calculated in the same way and we arrive at
\begin{align}
\left\{
\begin{array}{l}
\displaystyle
\dot{R_{1}}=-\frac{\alpha R_{2}(z_{1}-z_{2})}{\big( (R_{1}-R_{2})^{2}+(z_{1}-z_{2})^{2}\big)^{3/2}}, \\[7mm]
\displaystyle
\dot{z_{1}}=\frac{\beta}{R_{1}}+\frac{\alpha R_{2}(R_{1}-R_{2})}{\big( (R_{1}-R_{2})^{2}+(z_{1}-z_{2})^{2}\big)^{3/2}}, \\[7mm]
\displaystyle
\dot{R_{2}}=\frac{\alpha \beta R_{1}(z_{1}-z_{2})}{\big( (R_{1}-R_{2})^{2}+(z_{1}-z_{2})^{2}\big)^{3/2}}, \\[7mm]
\displaystyle
\dot{z_{2}}=\frac{1}{R_{2}}-\frac{\alpha \beta R_{1}(R_{1}-R_{2})}
{\big( (R_{1}-R_{2})^{2}+(z_{1}-z_{2})^{2}\big)^{3/2}}, \\[7mm]
(R_{1}(0),z_{1}(0),R_{2}(0),z_{2}(0))=(R_{1,0},z_{1,0},R_{2,0},z_{2,0}).
\end{array}\right.
\label{RZode}
\end{align}
First, we observe that \( z_{1}\) and \( z_{2}\) can be reduced to one variable, namely
\( W=z_{1}-z_{2}\). Furthermore, we see by direct calculation that 
\( \beta R_{1}^{2}+R_{2}^{2}\) is a conserved quantity. Hence, setting
\( d^{2}=\beta R_{1,0}^{2}+R_{2,0}^{2}\) with \( d>0\), 
we make the change of variables
\begin{align*}
R_{1}(t)=\frac{d}{\beta^{1/2}}
\cos( \theta (t)), \quad R_{2}(t)=d\sin( \theta(t))
\end{align*}
to further reduce the system. We then arrive at
\begin{align}
\left\{
\begin{array}{l}
\displaystyle
\dot{\theta} = 
\frac{\alpha \beta^{1/2}W}{
\big( \frac{d^{2}}{\beta}(\beta^{1/2}\sin \theta - \cos \theta)^{2}+W^{2}
 \big)^{3/2}} =:F_{1}(\theta ,W), \\[7mm]
\displaystyle
\dot{W}=\frac{\beta^{3/2}\sin \theta - \cos \theta}
{d\sin \theta \cos \theta}
-
\frac{
\alpha d^{2}(\sin \theta + \beta^{1/2}\cos \theta )
(\beta^{1/2}\sin \theta - \cos \theta)}
{\beta^{1/2}
\big( \frac{d^{2}}{\beta}(\beta^{1/2}\sin \theta - \cos \theta)^{2}+W^{2}
 \big)^{3/2}}
=: F_{2}(\theta ,W),
\end{array}\right.
\label{2dyn}
\end{align}
with initial data \( (\theta_{0},W_{0})\). Here,
\( W_{0}=z_{1,0}-z_{2,0}\) and \( \theta _{0}\) is determined uniquely from the relation
\begin{align*}
R_{1,0}=\frac{d}{\beta^{1/2}}\cos \theta_{0}, \quad R_{2,0}=d\sin \theta_{0}.
\end{align*}
Note that from our problem setting, \( (\theta_{0},W_{0}) \) 
is contained in 
the open set \( \Omega_{\beta} \subset \mathbf{R}^{2} \) given by
\begin{align*}
\Omega_{\beta}= \big\{ (\theta ,W)\in \mathbf{R}^{2} \ | \ 0<\theta <\frac{\pi}{2}, W\in \mathbf{R},
(\theta ,W)\neq (\theta_{\beta},0)\big\},
\end{align*}
where \( \theta_{\beta}\) is the unique solution of
\begin{align*}
\beta^{1/2}\sin \theta_{\beta} - \cos \theta_{\beta} =0,
\end{align*}
which is given explicitly by \( \theta_{\beta}=\arctan(1/\beta^{1/2}) \).
The excluded point in the above definition corresponds to the two filaments overlapping.
Since we can reconstruct the solution of (\ref{RZode}) from 
the solution \( (\theta(t),W(t))\) of (\ref{2dyn}) by 
\begin{align*}
R_{1}(t)=\frac{d}{\beta^{1/2}}\cos(\theta(t)), \quad R_{2}(t)=d\sin(\theta(t)),\\[3mm]
z_{1}(t)=\int^{t}_{0}\frac{\beta}{R_{1}(\tau)}
+
\frac{\alpha R_{2}(\tau )(R_{1}(\tau )-R_{2}(\tau))}
{ \big((R_{1}(\tau)-R_{2}(\tau))^{2}+W(\tau )^{2}\big)^{3/2}}
{\rm d}\tau , \\[3mm]
z_{2}(t)=\int^{t}_{0}\frac{1}{R_{2}(\tau)}-
\frac{\alpha \beta R_{2}(\tau )(R_{1}(\tau )-R_{2}(\tau))}
{ \big((R_{1}(\tau)-R_{2}(\tau))^{2}+W(\tau )^{2}\big)^{3/2}}
{\rm d}\tau,
\end{align*}
we focus on the solvability and behavior of the solution to system (\ref{2dyn}).
It can be checked by direct calculation that the system (\ref{2dyn}) is a
Hamiltonian system and the Hamiltonian \( {\mathcal H}\) is given by
\begin{align}
{\mathcal H}(\theta ,W)=\frac{1}{2d}\log\left(
\frac{(1-\sin \theta )^{\beta^{3/2}}(1- \cos \theta) }{(1+\sin \theta)^{\beta^{3/2}}(1+\cos \theta )}\right)
-
\frac{\alpha \beta^{1/2}}
{\big( \frac{d^{2}}{\beta}(\beta^{1/2}\sin \theta - \cos \theta)^{2}+W^{2}
 \big)^{1/2}}.
\label{hamil}
\end{align}
In other words, \( F_{1}=\frac{\partial {\cal H}}{\partial W}\) and 
\( F_{2}=-\frac{\partial {\cal H}}{\partial \theta }\).
Of course, the Hamiltonian is a conserved quantity of motion.
In this formulation, closed orbits revolving around the point
\( (\theta_{\beta},0)\) correspond to leapfrogging.
From here, we treat (\ref{2dyn}) as a two-dimensional dynamical system in 
\( \Omega_{\beta} \) with parameters
\( d\),\( \beta \), and \( \alpha \), and make use of many tools known for two-dimensional dynamical systems and 
Hamiltonian systems, for example in Hirsch and Smale \cite{37}, to 
determine the dynamics of the filaments.

\medskip
\medskip
We state our main theorems.
\begin{Th}
For any \( \alpha, d>0\), \( \beta \geq 1 \), and \( (\theta_{0},W_{0})\in \Omega_{\beta } \), 
there exists a unique time-global 
solution \( (\theta ,W)\in C^{1}\big( \mathbf{R}\big) 
\times C^{1}\big( \mathbf{R}\big)\) of {\rm (\ref{2dyn})}.
\label{TH1}
\end{Th}
\begin{Th}
In addition to the assumptions of Theorem {\rm \ref{TH1}}, if we assume \( 0<\alpha <1/3 \), then  
system {\rm (\ref{2dyn})} has two equilibrium points \( (\theta_{\ast},0) \) and \( (\theta_{\ast \ast},0)\)
with \( 0<\theta_{\ast}<\theta_{\beta}\) and \( \theta_{\beta}<\theta_{\ast \ast}<\pi/2 \), and 
the following two statements are equivalent.
\begin{description}
\item[\quad (i)] The solution with initial data
 \( (\theta_{0},W_{0}) \) is a 
leapfrogging solution. In other words, the solution curve is a closed orbit 
revolving around the point \( (\theta_{\beta},0)\).
\item[\quad (ii)] \( \theta_{0}\in (\theta_{\ast},\theta_{\ast \ast})\) and 
\( {\cal H}(\theta_{0},W_{0})<\min \{ {\cal H}(\theta_{\ast},0),{\cal H}(\theta _{\ast \ast},0) \}\).
\end{description}
\label{TH2}
\end{Th}
\begin{Rm}{\rm (Note on the assumption for \( \alpha \) in 
Theorem {\rm \ref{TH2}})}
Recall that \( \alpha >0 \) was given by
\( \alpha = 2\delta/\log(\frac{L}{\varepsilon}) \), where
\( \delta , \varepsilon >0\) were small parameters
with \( L>0\) fixed. These parameters were 
introduced in the course of the 
derivation of the model system {\rm (\ref{model})}. Hence, it is 
natural to assume that \( \alpha \) is small and also important that 
the smallness assumption for \( \alpha \) 
in Theorem {\rm \ref{TH2}} is 
independent of the parameters \( d\) and \( \beta \).
This allows us to treat different configurations of the filaments in the 
framework of one model, as opposed to models with different 
parameters depending on the configuration.
\end{Rm}

The rest of the section is devoted to the proof the above two theorems.

\medskip
\noindent
{\it Proof of Theorem {\rm \ref{TH1}}.}
Since \( F_{1}\) and \( F_{2}\) are smooth
in \( \Omega_{\beta}\), the time-local unique solvability is known. 
Suppose the maximum existence time \( T>0\) is finite. From the standard theory of dynamical systems,
for any compact set \( K\subset \Omega_{\beta}\), there exists \( t'\in [0,T)\) such that 
\( (\theta(t' ),W(t'))\not\in K\).
On the other hand, since the Hamiltonian is conserved, there exists \( \eta >0\) and 
\( r>0\) such that for all \( t\in [0,T)\),
\begin{align*}
(\theta(t),W(t))\in \big( [\eta, \frac{\pi}{2}-\eta ]\times \mathbf{R}\big)\setminus
B_{r}(\theta_{\beta},0),
\end{align*}
where \( B_{r}(\theta_{\beta},0)\) is the open ball with center \( (\theta_{\beta},0)\) and radius \( r\).
This follows from the fact that the Hamiltonian diverges to \( -\infty \) at 
\( \theta =0,\pi/2\) uniformly with respect to \( W\) and at the point \( (\theta_{\beta},0)\).
In particular, since the solution curve is uniformly away from the 
point \( (\theta_{\beta},0)\), there exists \( c_{0}>0\) such that
\begin{align*}
\frac{d^{2}}{\beta}(\beta^{1/2}\sin \theta(t) - \cos \theta(t))^{2}+W(t)^{2}\geq c_{0}
\end{align*}
for all \( t\in [0,T)\). Hence from the second equation in (\ref{2dyn}), we have
\begin{align*}
|\dot{W}| \leq \frac{\beta^{3/2}+1}{d\sin \eta \cos (\pi/2-\eta)}
+
\frac{\alpha d^2 (\beta^{1/2}+1)^{2}}{\beta^{1/2}c_{0}^{3/2}}=:M,
\end{align*}
which yields
\begin{align*}
|W(t)|\leq |W(0)|+ Mt \leq |W_{0}| + MT
\end{align*}
for all \( t\in [0,T)\). Finally, this shows that 
for all \( t\in [0,T)\), \( (\theta(t),W(t)) \) is contained in the compact set \( K' \) given by
\begin{align*}
K' = 
\big( [\eta, \frac{\pi}{2}-\eta ]\times [-|W_{0}|-MT,|W_{0}|+MT]\big)\setminus
B_{r}(\theta_{\beta},0),
\end{align*}
which is a contradiction. 
The same argument holds for \( t<0\) and hence, 
the solution exists globally in time and is defined for 
all \(t\in \mathbf{R}\).
\hfill \( \Box \)

\medskip

\medskip 

\noindent
{\it Proof of Theorem {\rm \ref{TH2}}.}
We divide the proof of Theorem \ref{TH2} into subsections.
First we prove that system (\ref{2dyn}) has exactly two equilibriums as stated in Theorem \ref{TH2}.
%
%

\subsection{Equilibriums of System (\ref{2dyn})}

From the form of \( F_{1}\), we see that an equilibrium can only exist on the line segment
\( (0,\pi/2)\times \{0\} \), and thus, we set
\( f(\theta ):=F_{2}(\theta ,0) \)
and investigate the zeroes of \( f\). 
First we consider the zeroes in the interval \( (0,\theta_{\beta})\).
Keeping in mind that \(\beta^{1/2}\sin \theta -\cos \theta <0\) in \( (0,\theta_{\beta})\), 
by a change of variable \( \theta = \arctan x\) we have
\begin{align*}
f(\arctan x) = 
\frac{(1+x^{2})^{1/2}g_{\alpha }(x)}{dx(\beta^{1/2}x-1)^{2}},
\end{align*}
where \( g_{\alpha }\) is given by
\begin{align*}
g_{\alpha}(x)=
\beta^{5/3}x^{3}-\beta(2\beta +1)x^{2}+\beta^{1/2}(\beta +2)x-1
+
\alpha \beta (x^{2}+\beta^{1/2}x)
\end{align*}
for \( x\in (0,1/\beta^{1/2})\).
We further make the change of variable \( y=\beta^{1/2}x\) for 
simplification and investigate the zeroes of the function \( h_{\alpha}\) given by
\begin{align*}
h_{\alpha}(y)=\beta y^{3}-(2\beta +1)y^{2}+(\beta +2)y - 1 +\alpha (y^{2}+\beta y)
\end{align*}
in the interval \( I_{1}=(0,1) \).
We treat \( h_{\alpha}\) as a perturbation of \( h_{0}\) given by
\begin{align*}
h_{0}(y)=\beta y^{3}-(2\beta +1)y^{2}+(\beta +2)y - 1,
\end{align*}
which is \( h_{\alpha }\) with \( \alpha =0\) and prove that 
\( h_{\alpha }\) has exactly one zero in \( I_{1}\).
We see from direct calculation that 
\( h_{0}\) has one local maximum and one local minimum at
\begin{align*}
y_{1}=\frac{\beta+2}{3\beta}, \quad y_{2}=1,
\end{align*}
respectively, and
\begin{align*}
h_{0}(y_{1})=\frac{4}{27\beta^{2}}(\beta-1)^{3}>0, \quad h_{0}(y_{2})=0.
\end{align*}
Since the zero at \( y_{2}\) is singular, we cannot directly apply the 
method of perturbation to \( h_{\alpha } \).
Instead, we analyze the positions of the local extrema for small \( \alpha >0\) to 
determine the number of zeroes of \( h_{\alpha }\).
First, we observe that the discriminant \( \Delta \) of the quadratic equation
\( h'_{\alpha}(y)=0 \) is given by
\begin{align*}
\Delta = 4[(1-3\alpha )\beta^{2}-2(1+2\alpha )\beta + (\alpha -1)]=:4\phi(\beta ).
\end{align*}
\( \phi(\beta )=0\) has two roots \( \beta_{\pm}\) given by
\begin{align*}
\beta_{-}=\frac{1+2\alpha -\sqrt{3\alpha (3-\alpha -\alpha^{2})}}{1-3\alpha },
\quad
\beta_{+}=
\frac{1+2\alpha + \sqrt{3\alpha (3-\alpha -\alpha^{2})}}{1-3\alpha }
\end{align*}
and under the assumption \( 0<\alpha <1/3 \), we see that
\begin{align*}
\phi(\beta) < 0 \ \mbox{for} \ 1\leq \beta < \beta_{+}, \quad 
\phi(\beta) \geq 0 \ \mbox{for} \ \beta_{+}\leq \beta,
\end{align*}
where we also used the fact that \( \phi(1)=-\alpha (9-\alpha )<0 \).
This shows that when \( 1\leq \beta<\beta_{+}\), \( \Delta <0\) which implies
\( h'_{\alpha }>0\) for \( y\in (0,1)\).
Since, \( h_{\alpha}(0)=-1\) and \( h_{\alpha }(1)=\alpha (1+\beta )>0\),
there is exactly one zero in \( I_{1}\).

When \( \beta_{+} \leq \beta\), the roots \( y_{\pm}\) of 
\( h'_{\alpha}(y)=0\) are given by
\begin{align*}
y_{\pm}=\frac{2\beta+1-\alpha \pm \sqrt{\phi(\beta)}}{3\beta},
\end{align*}
where \( y_{-}\) is the local maximum and \( y_{+}\) is the local minimum.
Since \( h_{\alpha}\) is a third order polynomial, it is 
sufficient to prove that \( h_{\alpha}(y_{+})>0\) to prove that 
\( h_{\alpha }\) has exactly one root.
We have 
\begin{align*}
y_{+}\geq \frac{1}{3\beta}(2\beta +1-\alpha )
&\geq \frac{1}{3\beta }(\beta +2+(\beta_{+}-1)-\alpha )\\[3mm]
&=
\frac{1}{3\beta }\big\{ \beta +2 + \frac{\alpha^{1/2}}{1-3\alpha}
[\big((3(3-\alpha -\alpha^{2})\big)^{1/2}
+5\alpha ^{1/2}-(1-3\alpha )\alpha^{1/2} ]\big\}\\[3mm]
&\geq \frac{\beta+2}{3\beta},
\end{align*}
which implies \( h_{0}(y_{+})\geq 0 \).
Finally, we have
\begin{align*}
h_{\alpha}(y_{+})=h_{0}(y_{+})+\alpha (y_{+}^{2}+\beta y_{+})>0
\end{align*}
which shows that \( h_{\alpha} \) also has exactly one root
when \( \beta_{+}\leq \beta \). Hence we have proven that 
for any \( \beta \geq 1 \) and \( 0<\alpha <1/3 \), 
\( h_{\alpha } \) has exactly one zero \( y_{\ast}\) in \( I_{1} \)
 and \( h'_{\alpha}(y_{\ast})>0\).
Hence, \( \theta_{\ast}=\arctan (y_{\ast}/\beta^{1/2}) \) is the 
desired zero of \( f(\theta) \) in the interval \( (0,\theta_{\beta})\) and we see that \( f'(\theta_{\ast})>0\).
By a similar argument, we see that there exists 
a unique \( \theta_{\ast \ast}\in (\theta_{\beta},\pi/2) \) such 
that \( f(\theta_{\ast \ast})=0\) and \( f'(\theta_{\ast \ast})>0\).
We note here that because \( \theta _{\ast}\) and 
\( \theta_{\ast \ast}\) are the only zeroes in the interval 
\( (0,\theta_{\beta})\) and \( (\theta_{\beta},\pi/2)\)
respectively, and \( f'(\theta_{\ast}),f'(\theta_{\ast \ast})>0\),
we have the following property for \( f(\theta )\).
\begin{align}
\begin{array}{ll}
f(\theta)<0, & \mbox{for} \ \theta \in (0,\theta_{\ast})\cup 
(\theta_{\ast \ast},\pi/2), \\[3mm]
f(\theta)>0, & \mbox{for} \ \theta \in (\theta_{\ast},\theta_{\beta})
\cup (\theta_{\beta},\theta_{\ast \ast}).
\end{array}
\label{propf}
\end{align}

%
%
%
%
%
%
%
%
\subsection{Analysis for Solutions with Initial Data of the Form 
\( (\theta_{0},0)\)}
Since a leapfrogging solution corresponds to a 
closed orbit revolving around the point \( (\theta_{\beta},0)\)
in \( \Omega_{\beta}\), 
a leapfrogging solution always crosses the lines
\( (0,\theta_{\beta})\times \{0\}\) and 
\( (\theta_{\beta},\pi/2)\times \{0\}\) in \( \Omega_{\beta}\).
To this end, we first characterize
the solutions with initial data of the form 
\( (\theta_{0},0)\), and prove that the condition given in Theorem \ref{TH2} is necessary and sufficient
for leapfrogging to occur.

First we prove that (ii) implies (i).
Set \( H_{\ast}:= \min \{{\cal H}(\theta_{\ast},0), 
{\cal H}(\theta_{\ast \ast},0)\} \).
Let \( \theta_{0}\in (\theta_{\ast},\theta_{\ast \ast})\) satisfy
\( {\cal H}(\theta_{0},0)<H_{\ast}\).
To make the situation more concrete, we further assume that 
\( {\cal H}(\theta_{\ast},0)>
{\cal H}(\theta_{\ast \ast },0)\) and make a 
remark on the case
\( {\cal H}(\theta_{\ast},0)\leq 
{\cal H}(\theta_{\ast \ast },0)\)
at the end.
From (\ref{propf}) and the fact that 
\( \frac{\partial {\cal H}}{\partial \theta}(\theta ,0)
=
-f(\theta ) \), we have
\begin{align}
\begin{array}{ll}
\displaystyle 
\frac{\partial {\cal H}}{\partial \theta}(\theta ,0)>0, 
& \mbox{for} \ \theta \in (0,\theta_{\ast})\cup 
(\theta_{\ast \ast},\pi/2), \\[5mm]
\displaystyle
\frac{\partial {\cal H}}{\partial \theta}(\theta ,0)<0, 
& \mbox{for} \ \theta \in (\theta_{\ast},\theta_{\beta})
\cup (\theta_{\beta},\theta_{\ast \ast}).
\end{array}
\label{propH}
\end{align}
Moreover, since \( {\cal H}(\theta_{\ast},0)>
{\cal H}(\theta_{\ast \ast},0)\), and 
\( {\cal H}(\theta ,0) \to -\infty \) monotonically as 
\( \theta \to \theta_{\beta}-\), there exists a 
unique \( \tilde{\theta} \in (\theta_{\ast},\theta_{\beta}) \)
such that 
\( {\cal H}(\tilde{\theta},0)=H_{\ast} \). 
This implies that \( \theta_{0}\in (\tilde{\theta},\theta_{\ast \ast})
\setminus \{\theta_{\beta }\}\).

%
%
%
%
%
%
%
%
We assume that \( \theta_{0}\in (\tilde{\theta},\theta_{\beta})\)
since the arguments for the case \( \theta_{0}\in (\theta_{\beta },
\theta_{\ast \ast})\) is the same.
We prove that the unique time-global solution
\( (\theta(t),W(t))\) starting from 
\( (\theta_{0},0)\) obtained in 
Theorem \ref{TH1}, which is defined for 
\( t\in \mathbf{R}\), is a closed orbit revolving around 
\( (\theta_{\beta},0) \).
First, we show that the solution is bounded. 
We observe that as a function of 
\( W \), the Hamiltonian achieves a minimum 
at \( W=0\) for each fixed \( \theta \).
Hence for all \( W\in \mathbf{R}\), we have
\begin{align*}
{\cal H}(\tilde{\theta},W)&\geq 
{\cal H}(\tilde{\theta},0) = H_{\ast}
>{\cal H}(\theta_{0},0), \\[3mm]
{\cal H}(\theta_{\ast \ast},W)&\geq
{\cal H}(\theta_{\ast \ast},0)=H_{\ast}
>{\cal H}(\theta_{0},0). 
\end{align*}
The above and from the conservation and 
continuity of the Hamiltonian, there
exists \( \eta >0\) and \( r>0\) such that
\begin{align*}
(\theta(t),W(t)) \in 
\big( 
[\tilde{\theta}+\eta , \theta_{\ast \ast}-\eta]\times 
\mathbf{R}\big)
\setminus
B_{r}(\theta_{\beta},0),
\end{align*}
for all \( t\in \mathbf{R} \).
Furthermore, if we set 
\begin{align*}
\phi(\theta ):= 
\frac{1}{2d}\log\left(
\frac{(1-\sin \theta )^{\beta^{3/2}}(1- \cos \theta) }{(1+\sin \theta)^{\beta^{3/2}}(1+\cos \theta )}\right),
\end{align*}
we see that as a function of \( \theta \), 
\( {\cal H}(\theta , W) \) converges to 
\( \phi \) uniformly as \( W\to \infty \).
Since we have 
\begin{align*}
\phi '(\theta )=
-\frac{(\beta^{3/2}\sin \theta -\cos \theta )}{d \cos \theta 
\sin \theta },
\end{align*}
we see that \( \phi \) achieves a maximum 
at \( \theta = \arctan(1/\beta^{3/2})=:\theta_{c}\) with
\( 0<\theta_{c}<\theta_{\beta}\) and 
\( \phi \) is monotone in the intervals \( (0,\theta_{c})\) and 
\( (\theta_{c},\pi/2)\).
If \( 0<\theta_{c}\leq \tilde{\theta}\), for 
\( \varepsilon_{1}>0\) given by 
\begin{align*}
\varepsilon_{1}=
\frac{\alpha \beta^{1/2}}
{
2\big\{\frac{d^{2}}{\beta}\big( \beta^{1/2}
\sin \theta_{\ast \ast}
-\cos \theta_{\ast \ast}\big)^{2} \big\}^{1/2}},
\end{align*}
there exists \( W_{1}>0\) such that for all 
\( \theta \in (\tilde{\theta},\theta_{\ast \ast} ) \),
and \( W> W_{1}\)
we have
\begin{align*}
{\cal H}(\theta ,W)>\phi (\theta )-\varepsilon_{1}
>\phi (\theta_{\ast \ast} )-2\varepsilon_{1}
={\cal H}(\theta_{\ast \ast},0) = H_{\ast}
>{\cal H}(\theta_{0},0).
\end{align*}
If \( \tilde{\theta}<\theta_{c}<\theta_{\beta} \),
choose \( \theta'\in \{ \tilde{\theta},\theta_{\ast \ast}\}\) 
so that \( \phi (\theta') = \min
\{ \phi (\tilde{\theta}),\phi(\theta_{\ast \ast})\}\).
Then for \( \varepsilon_{2}>0\) given by
\begin{align*}
\varepsilon_{2}=
\frac{\alpha \beta^{1/2}}
{
2\big\{\frac{d^{2}}{\beta}\big( \beta^{1/2}
\sin \theta'-\cos \theta'\big)^{2} \big\}^{1/2}},
\end{align*}
there exists \( W_{2}>0\) such that 
for all \( \theta \in (\tilde{\theta},\theta_{\ast \ast}) \) and 
\( W>W_{2}\), we have
\begin{align*}
{\cal H}(\theta ,W)>\phi (\theta)-\varepsilon_{2}
>\phi(\theta')-2\varepsilon_{2}
= {\cal H}(\theta',0) = H_{\ast}
> {\cal H}(\theta_{0},0).
\end{align*}
In either case, we see that the value of the Hamiltonian 
on the segment 
\( [\tilde{\theta},\theta_{\ast \ast}]\times 
\{W_{\ast}\}\), where \( W_{\ast}=\max \{ W_{1},W_{2}\} \), 
is strictly greater than \( {\cal H}(\theta_{0},0)\) and 
hence the solution curve cannot cross this segment.
Since the Hamiltonian is symmetric with respect to 
\( W=0\), we finally see that 
\begin{align*}
(\theta(t),W(t))\in 
\big( 
[\tilde{\theta}+\eta , \theta_{\ast \ast}-\eta]\times 
[-W_{\ast},W_{\ast}]\big)
\setminus
B_{r}(\theta_{\beta},0)=: K_{\ast},
\end{align*}
for all \( t\in \mathbf{R}\), and in particular, the solution is 
bounded. 

Next we set 
\begin{align*}
L_{0}:=\{(\theta,W)\in \Omega_{\beta} \ | \ 
{\cal H}(\theta,W)={\cal H}(\theta_{0},0) \}
\cap K_{\ast}.
\end{align*}
As a closed subset of the compact set \( K_{\ast}\), 
\( L_{0}\) is a compact subset of \( \Omega_{\beta}\).
From the conservation of the Hamiltonian
and the way we chose \( \eta \), \( r\), and
\( W_{\ast}\), we see that 
\( L_{0}\) is also an invariant set
and hence we have
\begin{align*}
L_{\omega}(\theta_{0},0)\subset L_{0},
\end{align*}
where \( L_{\omega}(\theta_{0},0)\) is the \( \omega \)-limit set 
of \( (\theta_{0},0)\). Since \( (\theta(t),W(t))\) is 
bounded for \( t>0\), it converges along some series
\( \{t_{n}\}^{\infty}_{n=1}\) with \( t_{n}\to \infty\)
as \( n\to \infty \),
and in particular, \( L_{\omega}(\theta_{0},0)\) is not empty. 
Since \( L_{\omega}(\theta_{0},0)\) is a non-empty compact set
and contains no equilibriums (recall that the equilibriums
\( (\theta_{\ast},0)\) and \( (\theta_{\ast \ast},0) \)
are outside the set \( L_{0}\)), it is a 
closed orbit by the Poincar\'e--Bendixson Theorem.
Moreover, the point \( (\theta_{\beta},0)\) is in the 
interior of this closed orbit, because if it is not,
then the closed orbit would enclose an open subset of 
\( \Omega_{\beta}\) in which an equilibrium must exist,
which leads to a contradiction.
This proves that \( L_{\omega}(\theta_{0},0)\) is a 
closed orbit revolving around \( (\theta_{\beta},0)\).
Since \( L_{\omega}(\theta_{0},0)\subset L_{0}\), there 
exists \( \theta_{1}\in (\tilde{\theta}+\eta , \theta_{\beta}) \)
and \( \theta_{2}\in (\theta_{\beta},\theta_{\ast \ast}-\eta) \)
such that \( (\theta_{1},0),(\theta_{2},0)\in 
L_{\omega}(\theta_{0},0)\).
The values \( \theta_{1}\) and \( \theta_{2}\) satisfying 
this property are unique in their respective intervals because 
the Hamiltonian is monotone along the line segments
\( [\tilde{\theta}+\eta , \theta_{\beta}]\times \{0\} \)
and \( [ \theta_{\beta},\theta_{\ast \ast}-\eta ]\times \{0 \}\).
This uniqueness implies that
\( \theta_{1}=\theta_{0}\), which proves that 
\( L_{\omega}(\theta_{0},0)\) coincides with the 
orbit starting from \( (\theta_{0},0)\).

In summary, we have proven that the orbit starting from 
\( (\theta_{0},0)\) is a closed orbit revolving around 
\( (\theta_{\beta},0)\) corresponding to a leapfrogging 
solution.
We further have the characterization
\begin{align*}
L_{\omega}(\theta_{0},0)=L_{0},
\end{align*}
which we prove by contradiction.
Suppose there exists \( (\overline{\theta},\overline{W})\in L_{0}\)
such that \( (\overline{\theta},\overline{W})\not\in 
L_{\omega}(\theta_{0},0)\).
We first see that \( \overline{W}\neq 0\), since 
\( (\overline{\theta},0)\in L_{0}\) implies
\( \overline{\theta}=\theta_{1} \ \mbox{or} \ \theta_{2}\),
which contradicts \( (\overline{\theta},0)\not\in
L_{\omega}(\theta_{0},0)\). 
Henceforth, we assume \( \overline{W}>0\) since the proof 
for the other case is the same. 
Now, if \( \overline{\theta}\in [\tilde{\theta}+\eta , \theta_{1}]\),
we have
\begin{align*}
{\cal H}(\overline{\theta},\overline{W})>
{\cal H}(\overline{\theta},0)\geq
{\cal H}(\theta_{1},0)
=
{\cal H}(\theta_{0},0)
\end{align*}
from the monotonicity of \( {\cal H}\) along the line
\( \{ \overline{\theta}\}\times \mathbf{R}\)
and the monotonicity along the line segment
\( [\overline{\theta},\theta_{1}]\times \{0\} \),
and this contradicts \( (\overline{\theta},\overline{W})\in L_{0}\).
The case \( \overline{\theta}\in 
[\theta_{2},\theta_{\ast \ast}-\eta ]\) leads to a contradiction
by the same argument.
If \( \overline{\theta}\in (\theta_{1},\theta_{2})\) and  
\( (\overline{\theta},\overline{W})\) is in the interior of the 
closed orbit \( L_{\omega}(\theta_{0},0)\), 
there exists \( \tilde{W}>\overline{W}\) such that
\( (\overline{\theta},\tilde{W})\in 
L_{\omega}(\theta_{0},0)\).
Then we have
\begin{align*}
{\cal H}(\overline{\theta},\overline{W})
<{\cal H}(\overline{\theta},\tilde{W})
=
{\cal H}(\theta_{0},0),
\end{align*}
which contradicts \( (\overline{\theta},\overline{W})\in L_{0}\).
Similarly, if \( (\overline{\theta},\overline{W})\) is 
outside of the closed orbit, there exists 
\( \tilde{W}<\overline{W}\) such that 
\( ( \overline{\theta},\tilde{W})\in L_{\omega}(\theta_{0},0)\).
Again, this implies the estimate
\begin{align*}
{\cal H}(\overline{\theta},\overline{W})
>{\cal H}(\overline{\theta},\tilde{W})
=
{\cal H}(\theta_{0},0),
\end{align*}
which contradicts \( (\overline{\theta},\overline{W})\in L_{0}\).
Hence we have \( L_{\omega}(\theta_{0},0)=L_{0}\).
We can express \( L_{0}\) as 
\begin{align*}
L_{0}=\{(\theta,W)\in \Omega_{\beta} \ | \ 
{\cal H}(\theta ,W)={\cal H}(\theta_{0},0)\}\cap
M,
\end{align*}
with \( M=[\theta_{\ast},\theta_{\ast \ast}]\times
\mathbf{R}\), because the value of the 
Hamiltonian on \( M\setminus K_{\ast}\) is different 
from \( {\cal H}(\theta_{0},0)\), and thus, 
replacing \( K_{\ast}\) with \( M\) does not add any points.
This expression will be utilized to derive the 
necessary and  sufficient condition for leapfrogging to 
occur for solutions with general initial data.

Finally, we make some remarks
on the case \( {\cal H}(\theta_{\ast},0) \not >
{\cal H}(\theta_{\ast \ast},0)\).
When 
\( {\cal H}(\theta_{\ast},0)={\cal H}(\theta_{\ast \ast},0)\),
the same proof holds with \( \tilde{\theta}=\theta_{\ast}\).
When
\( {\cal H}(\theta_{\ast},0)< {\cal H}(\theta_{\ast \ast},0)\),
there is a unique \( \hat{\theta}\in (\theta_{\beta},
\theta_{\ast \ast}) \) such that 
\( {\cal H}(\hat{\theta},0) = H_{\ast}\).
This \( \hat{\theta}\) plays the same role as 
\( \tilde{\theta}\), and the 
same arguments for the case 
\( {\cal H}(\theta_{\ast},0)< {\cal H}(\theta_{\ast \ast},0)\)
holds.

Next we prove that (i) implies (ii). 
Suppose that a solution starting from \( (\theta_{0},0)\) is a leapfrogging solution.
Since \( {\cal H}(\theta_{\ast},0)\) and \( {\cal H}(\theta_{\ast \ast},0) \) are the maximum
value of \( {\cal H}(\theta,0)\) in their respective intervals
\( (0,\theta_{\beta}) \) and \( (\theta_{\beta},\pi/2)\), in order for a solution curve to 
cross over the segments \( (0,\theta_{\beta})\times \{0\}\) and \( (\theta_{\beta},\pi/2)\times \{0\}\),
the value of the Hamiltonian on this solution curve must be less than or equal to the 
smaller of the two. In other words, \( {\cal H}(\theta_{0},0) \leq H_{\ast}\) holds.
If \( {\cal H}(\theta_{0},0)=H_{\ast}\) holds, the only possible points at which the solution curve can 
cross the segments \( (0,\theta_{\beta})\times \{0\}\) and \( (\theta_{\beta},\pi/2)\times \{0\}\)
are at the equilibrium points. This would result in the solution converging to one of the equilibrium points, 
and is not a leapfrogging solution.
Hence, for a leapfrogging solution, \( {\cal H}(\theta_{0},0)<H_{\ast} \) holds.

Furthermore, \( (\theta_{0},0)\) is not on the lines 
\( \{\theta_{\ast}\}\times\mathbf{R}\) or \( \{\theta_{\ast \ast}\}\times \mathbf{R}\) since the 
value of the Hamiltonian is greater than or equal to \( H_{\ast}\) along these lines.
Consequently, if \( \theta_{0}\in (0,\theta_{\ast})\cup (\theta_{\ast \ast}, \pi/2)\), the solution curve 
cannot cross over from one side of these lines to the other, which means that the solution is not a 
leapfrogging solution. This implies that \( \theta_{0}\in (\theta_{\ast},\theta_{\ast \ast})\), 
and condition (ii) holds.

\medskip
We summarize the conclusions of this 
subsection in the following lemma.
\begin{lm}
For initial data of the form \( (\theta_{0},0)
\in \Omega_{\beta}\), we have the 
following.
\begin{description}
\item[\quad (i)]If \( \theta_{0}\in 
(\theta_{\ast},\theta_{\ast \ast}) \) and 
\( {\cal H}(\theta_{0},0)<H_{\ast}\), then the solution starting 
from \( (\theta_{0},0)\) is a leapfrogging solution. Moreover,
the closed orbit \( L_{\omega}(\theta_{0},0)\) can be expressed as
\begin{align*}
L_{\omega}(\theta_{0},0)=
\{(\theta,W)\in \Omega_{\beta} \ | \ 
{\cal H}(\theta ,W)={\cal H}(\theta_{0},0)\}\cap
M,
\end{align*}
where \( M=[\theta_{\ast},\theta_{\ast \ast}]\times 
\mathbf{R}\).
\item[\quad (ii)]Otherwise, the solution is not a 
leapfrogging solution.
\end{description}

\label{lm}
\end{lm}
%
%

%
%
%
%
\subsection{Remarks on Solutions with General Initial Data}
Let \( (\theta_{0},W_{0})\in \Omega_{\beta} \) satisfy
\( \theta_{0}\in (\theta_{\ast},\theta_{\ast \ast})\) and 
\( {\cal H}(\theta_{0},W_{0})<H_{\ast} \).
Since \( {\cal H}(\theta,0) \) takes all values
between \( -\infty \) and \(H_{\ast} \) on the 
set \( (\theta_{\ast},\theta_{\beta})\cup
(\theta_{\beta},\theta_{\ast \ast})\), there exists 
\( \theta_{LF}\in (\theta_{\ast},\theta_{\beta})\cup
(\theta_{\beta},\theta_{\ast \ast}) \) such that 
\( {\cal H}(\theta_{LF},0) = {\cal H}(\theta_{0},W_{0}) \).
Moreover, from Lemma \ref{lm}, the orbit containing 
\( (\theta_{LF},0)\) is a closed orbit corresponding to 
a leapfrogging solution. Since
\begin{align*}
(\theta_{0},W_{0})\in 
\{(\theta,W)\in \Omega_{\beta} \ | \ 
{\cal H}(\theta ,W)={\cal H}(\theta_{LF},0)\}\cap
M,
\end{align*}
Lemma \ref{lm} implies that 
\( (\theta_{0},W_{0})\) is on the closed orbit 
containing \( (\theta_{LF},0)\) and hence,
the solution starting from \( (\theta_{0},W_{0})\) is a 
leapfrogging solution.

On the other hand, suppose either 
\( {\cal H}(\theta_{0},W_{0})\geq H_{\ast}\) or
\( \theta_{0} \not\in (\theta_{\ast},\theta_{\ast \ast})\) holds.
We prove that solution curves starting from these initial data
are not leapfrogging solutions.
If \( {\cal H}(\theta_{0},W_{0})\geq H_{\ast} \),
then the solution 
starting from \( ( \theta_{0},W_{0})\)
is not a leapfrogging solution since 
the value of the Hamiltonian of a leapfrogging solution is
strictly less than \( H_{\ast} \) from Lemma \ref{lm}.
If \( \theta_{0}\not\in (\theta_{\ast},\theta_{\ast \ast})\) holds,
we only need to consider the case
when \( {\cal H}(\theta_{0},W_{0})<H_{\ast}\) also 
holds.
Since \( {\cal H}(\theta_{0},W_{0})<H_{\ast} \),
\( \theta_{0}\in (0,\theta_{\ast})\cup 
(\theta_{\ast \ast},\pi/2)\) because the value of the Hamiltonian 
on the lines \( \{ \theta_{\ast}\}\times\mathbf{R}\) and 
\( \{ \theta_{\ast \ast}\}\times \mathbf{R}\) are greater than or 
equal to \( H_{\ast}\). Furthermore, since the Hamiltonian is 
conserved, the solution curve starting from 
\( (\theta_{0},W_{0})\) cannot cross over from one side of 
these lines to the other and hence, the solution is not a
leapfrogging solution.
This finishes the proof of Theorem \ref{TH2}.
\hfill \( \Box \)


\section{Leapfrogging for a Pair of Filaments with Vorticity Strengths of
Opposite Signs}
\setcounter{equation}{0}
We consider the case when the two filaments are a pair of coaxial circles with 
vorticity strengths of opposite signs. This amounts to considering
system (\ref{RZode}) with \( \beta <0 \). Again, since
the case \( -1<\beta <0\) is reduced to the case \( \beta \leq -1\) by 
renaming the filaments and rescaling the time variable, we assume 
\( \beta \leq -1 \) without loss of generality.
Setting \( \gamma = -\beta \), system (\ref{RZode}) reads
\begin{align}
\left\{
\begin{array}{l}
\displaystyle
\dot{R_{1}}=-\frac{\alpha R_{2}(z_{1}-z_{2})}{\big( (R_{1}-R_{2})^{2}+(z_{1}-z_{2})^{2}\big)^{3/2}}, \\[7mm]
\displaystyle
\dot{z_{1}}=-\frac{\gamma}{R_{1}}+\frac{\alpha R_{2}(R_{1}-R_{2})}{\big( (R_{1}-R_{2})^{2}+(z_{1}-z_{2})^{2}\big)^{3/2}}, \\[7mm]
\displaystyle
\dot{R_{2}}=-\frac{\alpha \gamma R_{1}(z_{1}-z_{2})}{\big( (R_{1}-R_{2})^{2}+(z_{1}-z_{2})^{2}\big)^{3/2}}, \\[7mm]
\displaystyle
\dot{z_{2}}=\frac{1}{R_{2}}+\frac{\alpha \gamma R_{1}(R_{1}-R_{2})}
{\big( (R_{1}-R_{2})^{2}+(z_{1}-z_{2})^{2}\big)^{3/2}}, \\[7mm]
(R_{1}(0),z_{1}(0),R_{2}(0),z_{2}(0))=(R_{1,0},z_{1,0},R_{2,0},z_{2,0}),
\end{array}\right.
\label{nRZode}
\end{align}
with \( \gamma \geq 1 \). 
Note that by the nature of the leapfrogging motion, 
if a solution of (\ref{nRZode}) corresponds to 
leapfrogging, \( R_{1}-R_{2}\) must change signs 
in a time-periodic pattern.
We see from direct calculation that
\( \gamma R_{1}^{2}-R_{2}^{2}\) is a conserved quantity.
This shows that \( (R_{1},R_{2})\) lies on the set defined by 
\( \gamma R_{1}^{2}-R_{2}^{2}=d \), where \( d=
\gamma R_{1,0}^{2}-R_{2,0}^{2}\), in the first quadrant of the 
\( R_{1}\)-\(R_{2}\) plane.
When \( d<0\), this is a hyperbola which approaches the 
line \( R_{2}=\gamma^{1/2}R_{1} \) from above and since \( \gamma \geq 1 \), 
\( R_{2}>R_{1}\) on the hyperbola. Hence, a solution cannot correspond to 
leapfrogging in this case.
If \( d=0\), \( (R_{1},R_{2})\) lies on the line 
\( R_{2}=\gamma^{1/2}R_{1} \).
When \( \gamma>1\), \( R_{2}>R_{1}\) and cannot correspond to leapfrogging.
When \( \gamma =1\), \( R_{1}=R_{2}\) throughout the motion, and 
the two filaments approaching would result in the collision of the two filaments.
This in itself is an interesting phenomenon, but is not a leapfrogging solution.
When \( d>0\), \( (R_{1},R_{2})\) lies on a hyperbola which approaches the line 
\( R_{2}=\gamma^{1/2}R_{1}\) from below.
If \( \gamma =1\), then \( R_{1}>R_{2}\) on this hyperbola and leapfrogging cannot occur.
If \( \gamma >1\), 
there is a possibility that the solution of (\ref{nRZode}) is a leapfrogging solution,
and we investigate in more detail.

First, we make the following change of variables.
\begin{align*}
R_{1}(t)=\left(\frac{d}{\gamma}\right)^{1/2}\cosh(\theta(t)), \quad 
R_{2}(t)=d^{1/2}\sinh(\theta(t)), \quad W(t)=z_{1}(t)-z_{2}(t),
\end{align*}
which yields
\begin{align}
\left\{
\begin{array}{l}
\displaystyle
\dot{\theta} = 
-\frac{\alpha \gamma^{1/2}W}{
\big( \frac{d}{\gamma}(\cosh \theta -\gamma^{1/2}\sinh \theta)^{2}+W^{2}
 \big)^{3/2}} =:G_{1}(\theta ,W), \\[7mm]
\displaystyle
\dot{W}=
-\frac{1}{d^{1/2}}
\left(
\frac{\gamma^{3/2}}{\cosh \theta}+\frac{1}{\sinh \theta}
\right)
+
\frac{
\alpha d(\sinh \theta - \gamma^{1/2}\cosh \theta )
(\cosh \theta - \gamma^{1/2}\sinh \theta)}
{\gamma^{1/2}
\big( \frac{d}{\gamma}(\cosh \theta -\gamma^{1/2}\sinh \theta)^{2}+W^{2}
 \big)^{3/2}} \\[7mm]
\phantom{W}=: G_{2}(\theta ,W),
\end{array}\right.
\label{n2dyn}
\end{align}
with initial data \( (\theta_{0},W_{0})\), which is determined by
\begin{align*}
R_{1,0}=\left(\frac{d}{\gamma}\right)^{1/2}\cosh(\theta_{0}), \quad
R_{2,0}=d^{1/2}\sinh(\theta_{0}), \quad 
W_{0}=z_{1,0}-z_{2,0}.
\end{align*}
The phase space \( \Omega_{\gamma } \subset \mathbf{R}^{2}\) is given by
\begin{align*}
\Omega_{\gamma}=\big\{ (\theta, W)\in \mathbf{R}^{2} \ | \ 
0<\theta <\infty, W\in \mathbf{R},(\theta,W)\neq (\theta_{\gamma },0)\big\},
\end{align*}
where \( \theta_{\gamma }\in (0,\infty) \) is the unique solution
of
\begin{align*}
\cosh \theta_{\gamma} -\gamma^{1/2}\sinh \theta_{\gamma} =0,
\end{align*}
given explicitly by \( \theta_{\gamma}=\artanh (1/\gamma^{1/2})\).
The point \( (\theta_{\gamma },0)\) corresponds to 
the two filaments overlapping.
In this formulation, leapfrogging solutions correspond to 
closed orbits revolving around the point
\( (\theta_{\gamma },0)\).
System (\ref{n2dyn}) is of Hamiltonian form and the Hamiltonian
\( {\cal G} \) is given by
\begin{align*}
{\cal G}(\theta,W) &=
\frac{1}{d^{1/2}}
\bigg(
2\gamma^{3/2}\arctan \big( \tanh(\theta/2)\big)+
\log(\tanh(\theta/2))
\bigg)\\[3mm]
& \hspace*{2cm}+
\frac{\alpha \gamma^{1/2}}
{\big( \frac{d}{\gamma}(\cosh \theta -\gamma^{1/2}\sinh \theta)^{2}+W^{2} \big)^{1/2}}.
\end{align*}

\medskip
We state our main theorems.
\begin{Th}
For any \( d>0\), \( \alpha >0\), \( \gamma > 1 \), and 
initial data \( (\theta_{0},W_{0})\in \Omega_{\gamma}\), system 
{\rm (\ref{n2dyn})} has a unique time-global solution
\( (\theta ,W)\in C^{1}(\mathbf{R})\times C^{1}(\mathbf{R}) \).
\label{TH3}
\end{Th}
\begin{Th}
Let \( 0<\alpha <1/3 \).
There exists \( \gamma_{\ast}\in (1,\infty)\) such that for any \( d>0\), the following holds.
When \( 1 < \gamma \leq \gamma_{\ast}\), system {\rm (\ref{n2dyn})} has no equilibriums and the following two statements are 
equivalent.
\begin{description}
\item[\quad (i)] The solution starting from \( (\theta_{0},W_{0})\in \Omega_{\gamma}\) is a leapfrogging solution.
\item[\quad (ii)] \( \displaystyle  {\cal G}(\theta_{0},W_{0})>\frac{\pi \gamma^{3/2}}{2d^{1/2}} \).
\end{description}
When \( \gamma_{\ast}<\gamma \), system {\rm (\ref{n2dyn})} has one equilibrium 
\( (\theta_{\ast},0)\) with \( \theta_{\ast}\in (\theta_{\gamma},\infty)\), and the following two statements are 
equivalent.
\begin{description}
\item[\quad (iii)] The solution starting from \( (\theta_{0},W_{0})\in \Omega_{\gamma}\) is a leapfrogging solution.
\item[\quad (iv)] \( \theta_{0}\in (\overline{\theta},\theta_{\ast})\) and \( {\cal G}(\theta_{0},W_{0})>
{\cal G}(\theta_{\ast},0)\).
\end{description}
Here, \( \overline{\theta}\in (0,\theta_{\gamma})\) is the unique value satisfying
\( {\cal G}(\overline{\theta},0)={\cal G}(\theta_{\ast},0)\).
\label{TH4}
\end{Th}
\noindent
{\it Proof of Theorem {\rm \ref{TH3}}.} 
The arguments of the proof are the same as that of Theorem \ref{TH1} and hence, we only give the details for the essential parts.
Let \( (\theta_{0},W_{0})\in \Omega_{\gamma }\) and set \( G_{0}:= {\cal G}(\theta_{0},W_{0})\). 
Since the unique existence of the time-local solution is known, we give an a priori estimate of the solution 
\( ( \theta(t),W(t))\)
on the time interval \( [0,T)\) to show that the solution is global in time.
We first see that since \( \cosh \theta -\gamma^{1/2}\sinh \theta >0\) in \( (0,\theta_{\gamma})\),
\begin{align*}
G_{2}(\theta,0)=
-\frac{1}{d^{1/2}}
\left\{
\frac{\gamma^{3/2}\sinh \theta + \cosh \theta}{\sinh \theta \cosh \theta}
-
\frac{\alpha \gamma (\sinh \theta -\gamma^{1/2}\cosh \theta)}
{(\cosh \theta-\gamma^{1/2}\sinh \theta)^{2}}
\right\}
<0,
\end{align*}
which implies \( \frac{\partial {\cal G}}{\partial \theta}(\theta,0) >0\) in \( (0,\theta_{\gamma})\).
By direct calculation, we see that  
\( {\cal G}(\theta,0)\to -\infty \) as \( \theta \to 0\) and \( {\cal G}(\theta,0)\to \infty \) as 
\( \theta \to \theta_{\gamma }\). This along with the monotonicity of \( {\cal G}(\theta,0)\) implies that 
there exists a unique \( \theta_{1}\in (0,\theta_{\gamma })\) such that 
\( {\cal G}(\theta_{1},0) = G_{0}/2\). Furthermore, \( {\cal G}(\theta, W)\) is a
strictly decreasing function of \( |W|\) for each fixed \( \theta \). 
Hence, the value of \( {\cal G}\) is strictly less than \( G_{0}\) along the line 
\( \{\theta_{1}\}\times \mathbf{R}\). Then, the conservation of the Hamiltonian implies that 
\( \theta_{1}\leq \theta(t)\) for any \( t\in [0,T)\). Additionally, since 
\( {\cal G}(\theta,W) \to \infty \) as \( (\theta,W)\to (\theta_{\gamma},0)\), 
there exists \( r>0\) such that 
\begin{align*}
(\theta(t),W(t))\in \big( [\theta_{1},\infty)\times \mathbf{R}\big)\setminus B_{r}(\theta_{\gamma},0)
\end{align*}
for all \( t\in [0,T)\). This implies that for some \( c_{0}>0\),
\begin{align*}
\frac{d}{\gamma}\big( \cosh \theta(t)-\gamma^{1/2}\sinh \theta(t)\big)^{2}+W(t)^{2} \geq c_{0}
\end{align*}
for all \( t\in [0,T)\). Hence we have
\begin{align*}
|\dot{W}| &\leq \frac{1}{d^{1/2}}
\left(
\gamma^{3/2}+\frac{1}{\sinh \theta_{1}}
\right)
+
\frac{
\alpha d|\sinh \theta - \gamma^{1/2}\cosh \theta |
|\cosh \theta - \gamma^{1/2}\sinh \theta|}
{\gamma^{1/2}
\big( \frac{d}{\gamma}(\cosh \theta -\gamma^{1/2}\sinh \theta)^{2}+W^{2}
\big)^{3/2}}.
\end{align*}
Note that the second term is bounded regardless of the size of \( \theta\). Indeed, 
for \( 0 \leq \theta \leq M \) with \( \theta_{\gamma}<M\), we have
\begin{align*}
&\frac{
\alpha d|\sinh \theta - \gamma^{1/2}\cosh \theta |
|\cosh \theta - \gamma^{1/2}\sinh \theta|}
{\gamma^{1/2}
\big( \frac{d}{\gamma}(\cosh \theta -\gamma^{1/2}\sinh \theta)^{2}+W^{2}
\big)^{3/2}} \\[3mm]
&\qquad \leq 
\frac{\alpha d(\sinh M+\gamma^{1/2}\cosh M)(\cosh M+\gamma^{1/2}\sinh M)}
{\gamma^{1/2}c_{0}^{3/2}}
=:M_{1}.
\end{align*}
For \( M<\theta \), we have 
\begin{align*}
&\frac{
\alpha d|\sinh \theta - \gamma^{1/2}\cosh \theta |
|\cosh \theta - \gamma^{1/2}\sinh \theta|}
{\gamma^{1/2}
\big( \frac{d}{\gamma}(\cosh \theta -\gamma^{1/2}\sinh \theta)^{2}+W^{2}
\big)^{3/2}}\\[3mm]
& \qquad \leq 
\frac{
\alpha \gamma^{1/2}|\sinh \theta - \gamma^{1/2}\cosh \theta |
|\cosh \theta - \gamma^{1/2}\sinh \theta|}
{
d^{1/2} |\cosh \theta -\gamma^{1/2}\sinh \theta|^{3}
}\\[3mm]
& \qquad \leq 
\frac{2\alpha \gamma^{1/2}
\big|
(1-\gamma^{1/2})e^{\theta}-(1+\gamma^{1/2})e^{-\theta}
\big|
\big|
(1-\gamma^{1/2})e^{\theta}+(1+\gamma^{1/2})e^{-\theta}
\big|
}
{d^{1/2}
\big(
(\gamma^{1/2}-1)e^{\theta}-(1+\gamma^{1/2})e^{-\theta}
\big)^{3}
}\\[3mm]
& \qquad =
\frac{2\alpha \gamma^{1/2}e^{-\theta}
\big|
(1-\gamma^{1/2})-(1+\gamma^{1/2})e^{-2\theta}
\big|
\big|
(1-\gamma^{1/2})+(1+\gamma^{1/2})e^{-2\theta}
\big|
}
{
d^{1/2}
\big(
(\gamma^{1/2}-1)-(1+\gamma^{1/2})e^{-2\theta}
\big)^{3}
}\\[3mm]
& \qquad \leq
\frac{4\alpha \gamma^{1/2}e^{-M}
\big(
(\gamma^{1/2}-1)+(1+\gamma^{1/2})e^{-2M}
\big)^{2}
}
{
d^{1/2}
\big(
(\gamma^{1/2}-1)-(1+\gamma^{1/2})e^{-2M}
\big)^{3}
}
=:M_{2}.
\end{align*}
Hence we have
\begin{align*}
|\dot{W}| \leq 
\frac{1}{d^{1/2}}
\left(
\gamma^{3/2}+\frac{1}{\sinh \theta_{1}}
\right)
+
\max\{M_{1},M_{2}\}
=:M_{0},
\end{align*}
which yields
\begin{align*}
|W(t)| \leq |W_{0}|+MT
\end{align*}
for all \( t\in [0,T)\). This further yields
\begin{align*}
|\dot{\theta}| \leq 
\frac{\alpha \gamma^{1/2}
\big(
|W_{0}|+MT
\big)}
{
c_{0}^{3/2}
}
=:C_{0},
\end{align*}
which implies
\begin{align*}
|\theta(t)|\leq |\theta_{0}|+C_{0}T
\end{align*}
for all \( t\in [0,T)\), and these estimates for \( \theta(t)\) and \( W(t)\) are the desired a priori estimates.\\
\hfill \( \Box\)

\medskip

\noindent
{\it Proof of Theorem {\rm \ref{TH4}}}. 
We first prove the existence of \( \gamma_{\ast}\in (1,\infty)\) as stated in the theorem. 

%
%
%
%
%
%
\subsection{The Existence of \( \gamma_{\ast}\in (1,\infty) \)}

From the form of system (\ref{n2dyn}), an equilibrium can only exist on the line
\( W=0\), and hence we look for zeroes of \( g(\theta):=G_{2}(\theta,0)\).
We have already seen that \( g(\theta)<0\) in \( (0,\theta_{\gamma})\), so we 
consider \( g(\theta)\) in the interval \( (\theta_{\gamma },\infty)\).
By the change of variable \( \theta = \artanh x\), we have
\begin{align*}
g(\artanh x)=
-\frac{(1+x)^{1/2}(1-x)^{1/2}}{d^{1/2}x(1-\gamma^{1/2}x)^{2}}g_{\alpha}(x),
\end{align*}
where \( x\in (1/\gamma^{1/2} ,1) \) and \( g_{\alpha }\) is given by
\begin{align*}
g_{\alpha}(x)=\gamma^{5/2}x^3 +(1-2\gamma)\gamma x^2 + (\gamma-2)\gamma^{1/2}x+1+\alpha \gamma^{1/2}x(\gamma^{1/2}x-\gamma^{3/2}).
\end{align*}
We further set \( y=\gamma^{1/2}x\) to simplify and obtain
\begin{align*}
h_{\alpha}(y)=\gamma y^{3}+(1-2\gamma)y^{2}+(\gamma -2)y+1 + \alpha y(y-\gamma)
\end{align*}
for \( y\in (1,\gamma^{1/2})\).
We see by direct calculation that \( h_{\alpha}\) takes local maximum at \( y_{-}\) and 
local minimum at \( y_{+}\), each given by
\begin{align*}
y_{\pm} = 
\frac{-(1-2\gamma +\alpha)\pm \sqrt{(\gamma+1)^{2}+\alpha (3\gamma^{2}-4\gamma +2+\alpha )}}
{3\gamma}.
\end{align*}
Note that \( 3\gamma^{2}-4\gamma +2+\alpha>0\) for all \( \gamma >1\) and in particular,
\( y_{+}-y_{-}\geq 2/3\) for all \( \gamma>1 \) and \( \alpha >0\). 
We further have
\begin{align*}
y_{-}<\frac{-1+2\gamma -\alpha}{3\gamma} =\frac{2}{3}-\frac{(1+\alpha)}{3\gamma}<1.
\end{align*}
Hence, regardless of the exact value of \( y_{+}\), it is sufficient to 
evaluate the value of \( h_{\alpha}\) at \( y=1, \gamma^{1/2}\) to determine the 
behavior of \( h_{\alpha}\) in the interval \( (1,\gamma^{1/2})\).
We have
\begin{align*}
h_{\alpha}(1)=\alpha(1-\gamma) <0
\end{align*}
and
\begin{align*}
h_{\alpha}(\gamma^{1/2}) = (\gamma^{1/2}-1)^{2}(\gamma^{3/2}+1)+\alpha \gamma (1-\gamma^{1/2}).
\end{align*}
Setting \( \eta := \gamma^{1/2}\), we have
\begin{align*}
h_{\alpha}(\eta)=(\eta-1)\big( \eta^{4}-3\eta^3-\alpha \eta^2 + \eta -1\big)
=:(\eta-1)\phi(\eta).
\end{align*}
After some differentiations, we can conclude that \( \phi \) is monotonically increasing in \( \eta>1\),
and 
\begin{align*}
\phi(1)=-\alpha , \qquad \phi(\eta)\to \infty \ (\eta \to \infty),
\end{align*}
which implies that there exists a unique \( \eta_{\ast}\in (1,\infty)\) such that \( \phi(\eta_{\ast})=0\).
Setting \( \gamma_{\ast}=\eta_{\ast}^2\), we see that if 
\( 1<\gamma \leq \gamma_{\ast}\), \( h_{\alpha}(\gamma^{1/2})\leq 0\) and hence \( h_{\alpha}<0\) in 
\( (1,\gamma^{1/2})\). This implies that \( g_{\alpha}<0\) in \( (\theta_{\gamma},\infty)\) and hence,
\( g>0\) in \( (\theta_{\gamma},\infty)\). This shows that there are no equilibriums in this case.

If \( \gamma_{\ast}<\gamma \), \( h_{\alpha}(\gamma^{1/2})>0\) which implies that
there exists a unique \( y_{\ast}\in (1,\gamma^{1/2})\) such that
\( h_{\alpha}(y_{\ast})=0\), and \( h_{\alpha}<0\) in \( (1,y_{\ast})\) and \( h_{\alpha}>0\) in 
\( (y_{\ast},\gamma^{1/2})\).
Then, setting \( \theta_{\ast}=\artanh(y_{\ast}/\gamma^{1/2}) \), 
\( (\theta_{\ast},0)\) is the unique equilibrium.

\medskip
The above arguments give the following profile for the Hamiltonian.
When \( 1<\gamma \leq \gamma_{\ast}\),
\begin{align}
\left\{
\begin{array}{ll}
\displaystyle \frac{\partial {\cal G}}{\partial \theta}(\theta,0)>0, & \theta \in (0,\theta_{\gamma}), \\[5mm]
\displaystyle \frac{\partial {\cal G}}{\partial \theta}(\theta,0)<0, & \theta \in (\theta_{\gamma},\infty).
\end{array}\right.
\label{ghamil1}
\end{align}
When \( \gamma_{\ast}<\gamma\),
\begin{align}
\left\{
\begin{array}{ll}
\displaystyle \frac{\partial {\cal G}}{\partial \theta}(\theta,0)>0, & \theta \in (0,\theta_{\gamma})
\cup (\theta_{\ast},\infty), \\[5mm]
\displaystyle \frac{\partial {\cal G}}{\partial \theta}(\theta,0)<0, & \theta \in (\theta_{\gamma},\theta_{\ast}).
\end{array}\right.
\label{ghamil2}
\end{align}
We also see from direct calculation that
\begin{align}
{\cal G}(\theta,0)\to -\infty \ (\theta \to 0), \quad {\cal G}(\theta,0)\to \infty \ (\theta \to \theta_{\gamma}).
\label{ghamil3}
\end{align}
These properties will be used later.

\medskip

Now we divide the proof of Theorem \ref{TH4} into the cases \( 1<\gamma \leq \gamma_{\ast}\) and 
\( \gamma_{\ast}<\gamma \). Like in Section 3, we determine the behavior of the solutions
starting from an initial data of the form \( (\theta_{0},0)\).
The arguments for solutions starting from a general initial data are similar to that given in Section 3 and 
we will only give a brief remark on the matter at the end.
%
%
%
%
%
%
%
\subsection{The Case \( 1<\gamma \leq \gamma_{\ast}\)}
\label{secs}
We prove that condition (i) implies (ii) and vice versa. 
Suppose we have a leapfrogging solution, which in other words, is a solution for which 
the solution curve is a closed orbit
revolving around the point \( (\theta_{\gamma},0)\) in the phase space \( \Omega_{\gamma} \).
For this to happen, the orbit must cross over the line 
\( (\theta_{\gamma},\infty)\times \{0\}\).
Since
\begin{align*}
{\cal G}(\theta,0)\to \frac{\pi \gamma^{3/2}}{2d^{1/2}}
\end{align*}
as \( \theta \to \infty\), the second property in 
(\ref{ghamil1}) and the conservation of the Hamiltonian asserts that
\( {\cal G}(\theta_{0},0)>\frac{\pi \gamma^{3/2}}{2d^{1/2}}\) must hold in order for the 
solution curve to cross over the line \( (\theta_{\gamma},\infty)\times \{0\}\).

Now, suppose \( G_{0}:={\cal G}(\theta_{0},0)>\frac{\pi \gamma^{3/2}}{2d^{1/2}}\) and consider the solution 
\( (\theta(t),W(t))\) starting from \( (\theta_{0},0)\).
First there exists \( \theta_{1}\in (0,\theta_{\gamma})\) and \( \theta_{2}\in (\theta_{\gamma},\infty)\),
both unique in their respective intervals, such that
\begin{align*}
{\cal G}(\theta_{1},0)={\cal G}(\theta_{2},0)=G_{0}.
\end{align*}
Note that either \( \theta_{1}\) or \( \theta_{2}\) is \( \theta_{0}\).
Set \( \zeta:= \theta_{1}/2 \). From the conservation of the Hamiltonian, there exists \( r>0\) such that
\begin{align*}
(\theta(t),W(t))\in \big( [\theta_{1}-\zeta, \theta_{2}+\zeta]\times \mathbf{R}\big)\setminus B_{r}(\theta_{\gamma},0)
\end{align*}
for all \( t\in \mathbf{R}\).
Furthermore, if we set 
\begin{align*}
\psi (\theta) :=
\frac{1}{d^{1/2}}
\big\{
2\gamma^{3/2}\arctan\big( \tanh(\theta/2)\big) + \log\big(\tanh(\theta/2)\big)
\big\},
\end{align*}
we see that \( {\cal G}(\theta, W)\) converges to \( \psi (\theta)\) 
uniformly as \( |W|\to \infty \) and this convergence is monotonically decreasing. Since we also have
\( \psi '>0\) in \( (0,\infty)\), for 
\( \varepsilon_{1}=\frac{1}{2}(G_{0}-\psi(\theta_{2}+\zeta)) \) there exists \( W_{\ast}>0\) such that for
all \( \theta \in [\theta_{1}-\zeta, \theta_{2}+\zeta ]\) and \( W\) satisfying \( |W|>W_{\ast}\), we have
\begin{align*}
{\cal G}(\theta, W)<\psi(\theta)+\varepsilon_{1}=\psi (\theta)+\frac{1}{2}\big(G_{0}-\psi (\theta_{2}+\zeta)\big)
<G_{0}.
\end{align*}
Again, the conservation of the Hamiltonian implies that 
\begin{align*}
(\theta(t),W(t))\in \big( [\theta_{1}-\zeta,\theta_{2}+\zeta]\times 
[-W_{\ast},W_{\ast}]\big)\setminus B_{r}(\theta_{\gamma},0)=:K.
\end{align*}
From here, the proof is the same as that of Section 3. The set 
\( L_{1}\subset \Omega_{\gamma}\) given by
\begin{align*}
L_{1}:= \{ (\theta,W)\in \Omega_{\gamma} \ | \ {\cal G}(\theta,W)=G_{0}\}\cap K
\end{align*}
is a compact invariant set which shows that
\( L_{\omega}(\theta_{0},0)\), the \( \omega \)-limit set of \( (\theta_{0},0)\), 
satisfies \( L_{\omega}(\theta_{0},0) \subset L_{1}\). Moreover, since \( L_{\omega}(\theta_{0},0)\)
is non-empty, compact, and contains no equilibriums, it is a closed orbit enclosing the point
\( (\theta_{\gamma},0)\). It follows that \( L_{1}=L_{\omega}(\theta_{0},0)\) and 
\( L_{\omega}(\theta_{0},0)\) coincides with the solution curve 
starting from \( (\theta_{0},0)\). Hence, the solution starting from \( ( \theta_{0},0) \) is a leapfrogging solution.

%
%
%
%
%
%
%
\subsection{The Case \( \gamma_{\ast}<\gamma \)}
\label{secl}
We prove that (iii) implies (iv) and vice versa. 
From (\ref{ghamil2}), we see that 
\( {\cal G}(\theta,0)\) takes its minimum value at \( \theta=\theta_{\ast}\) in the interval 
\( (\theta_{\gamma}, \infty) \).
We also see from (\ref{ghamil2}) and (\ref{ghamil3}) that there is a unique \(\overline{\theta}\in (0,\theta_{\gamma})\) such that 
\begin{align*}
{\cal G}(\overline{\theta},0)={\cal G}(\theta_{\ast},0).
\end{align*}

Suppose the solution starting from \( (\theta_{0},0)\) is a leapfrogging solution.
Since a leapfrogging solution must cross over the line
\( (\theta_{\gamma},\infty)\times \{0\}\), the conservation of the Hamiltonian implies that 
\( {\cal G}(\theta_{0},0)>{\cal G}(\theta_{\ast},0)\). 
This further implies that \( \overline{\theta}<\theta_{0}\) because
\( {\cal G}(\theta, 0) <{\cal G}(\theta_{0},0)\) for any \( \theta \in (0,\overline{\theta})\).
If \( \theta_{\ast}<\theta_{0}\), the solution curve starting from \( (\theta_{0},0)\) cannot
cross over the line \( \{\theta_{\ast}\}\times \mathbf{R}\), and in turn is not a leapfrogging solution.
This shows that \( \theta_{0}<\theta_{\ast}\) must hold and hence condition (iv) is satisfied.

Now, suppose condition (iv) holds for \( (\theta_{0},0)\in \Omega_{\gamma}\), and 
consider the solution \( (\theta(t),W(t))\) starting from \( (\theta_{0},0)\).
The conservation of the Hamiltonian asserts that there exists \( r>0\) such that
\begin{align*}
(\theta(t),W(t))\in \big( [\overline{\theta}-\zeta,\theta_{\ast}+\zeta]\times \mathbf{R}\big) \setminus B_{r}(\theta_{\gamma},0)
\end{align*}
for all \( t\in \mathbf{R}\), where \( \zeta=\frac{1}{2}\min \{ \theta_{0}-\overline{\theta},\theta_{\ast}-\theta_{0}\}\).
Again, from the monotone and uniform convergence of \( {\cal G}(\theta,W)\) as \( |W|\to \infty\), we see that 
there exists \( W_{\ast}>0\) such that
\begin{align*}
(\theta(t),W(t))\in \big( [\overline{\theta}-\zeta,\theta_{\ast}+\zeta]\times [-W_{\ast},W_{\ast}]\big) 
\setminus B_{r}(\theta_{\gamma},0)
\end{align*}
for all \( t\in \mathbf{R}\). The rest is the same as the previous case 
and we have the conclusion that the \( \omega\)-limit set \( L_{\omega}(\theta_{0},0)\) of 
\( (\theta_{0},0)\) is the desired closed orbit corresponding to a leapfrogging solution.

%
%
%
%
%
%
\subsection{Remarks on Solutions with General Initial Data}
We give a brief remark for solutions starting from general initial data
\( (\theta_{0},W_{0})\in \Omega_{\gamma}\).
From the proof given so far, we see that in the case 
\( 1<\gamma \leq \gamma_{\ast}\), the collection of all closed orbits coincides with the collection of sets 
\( \{L(G)\}\) where the set \( L(G)\) is
given by
\begin{align*}
L(G):= \{ (\theta,W)\in \Omega_{\gamma} \ | \ {\cal G}(\theta,W)=G\}
\end{align*}
for \( G>\frac{\pi \gamma^{3/2}}{2d^{1/2}}\). For a solution starting from \( (\theta_{0},W_{0})\) to be a 
leapfrogging solution, it is necessary and sufficient that 
\( (\theta_{0},W_{0}) \in L(G)\) for some \( G>\frac{\pi \gamma^{3/2}}{2d^{1/2}} \).
Similarly, in the case 
\( \gamma_{\ast}<\gamma\), the collection of all closed orbits coincides with the collection of sets \( \{K(G)\} \), where
the set \( K(G)\) is given by 
\begin{align*}
K(G):= \big( [\overline{\theta},\theta_{\ast}]\times \mathbf{R} \big) \cap L(G)
\end{align*}
for \( G>{\cal G}(\theta_{\ast},0)\). Again, for the solution to be a leapfrogging solution,
it is necessary and sufficient that 
\( (\theta_{0},W_{0})\in K(G)\) for some \( G>{\cal G}(\theta_{\ast},0)\).
In both cases, the conditions given in Theorem \ref{TH4} is a reinterpretation 
of these facts. This finishes the proof of Theorem \ref{TH4}. \hfill \( \Box\)

%
%
%
%
%
%

\section{Discussions and Conclusions}
\setcounter{equation}{0}
We make some comparisons with existing models and make concluding remarks.
%
%
%
%
%
%
%
%
\subsection{Comparison with Existing Results}
We make a comparison between the 
results by Borisov, Kilin, and Mamev \cite{39}.
In \cite{39}, they consider the following model system.
\begin{align}
\left\{
\begin{array}{l}
\displaystyle \dot{R}_{i}=-\frac{1}{R_{i}}\frac{\partial}{\partial Z_{i}}
\sum_{j\neq i}\Gamma_{j}G(R_{i},Z_{i},R_{j},Z_{j})\\[7mm]
\displaystyle
\dot{Z_{i}}=\frac{1}{R_{i}}\frac{\partial}{\partial R_{i}}
\left(
\sum_{j\neq i}\Gamma_{j}G(R_{i},Z_{i},R_{j},Z_{j})
\right)
+
\frac{\Gamma_{i}}{4\pi R_{i}}
\left(
\log\frac{8R_{i}}{a_{i}}-\frac{3}{4}
\right),
\end{array}\right.
\label{ODE}
\end{align}
where \( i=1,2 \) is the index for the 
two rings, \( R_{i}\) are the radii of the 
rings, \( Z_{i}\) are the distances along the common axis of 
symmetry, \( \Gamma_{i}\) are the vorticity strengths of the rings,
\( a_{i}\) are the radii of the cross-section of the cores, 
which is taken to be a constant, and \( G\) is
given by
\begin{align*}
G(z,r,\tilde{z},\tilde{r}) =
\frac{(r\tilde{r})^{1/2}}{2\pi }
\left(
\big( \frac{2}{k}-k\big)K(k)-\frac{2}{k}E(k)\right),
\quad k=\left(\frac{4r\tilde{r}}{(z-\tilde{z})^{2}+
(r+\tilde{r})^{2}}\right)^{1/2}.
\end{align*}
\( K(k)\) and \( E(k)\) are the complete elliptic integrals of the first and second kind
given by
\begin{align*}
K(k)=\int^{1}_{0}\frac{{\rm d}x}{\sqrt{1-x^{2}}\sqrt{1-k^{2} x^{2}}}, \quad
E(k)=\int^{1}_{0}\frac{\sqrt{1-k^{2} x^{2}}}{\sqrt{1-x^{2}}} {\rm d}x.
\end{align*}
In \cite{39}, they analyze (\ref{ODE}) to determine all the possible types of
motion.
The model used in \cite{39}
is derived as a system of ordinary differential equations (ODE), in other words, 
they only consider thin vortex rings as opposed to filaments with general shapes.
The advantages of (\ref{ODE}) over our
model (\ref{model}) is that it can incorporate the effect of the change in shape
of the core and the effect of the change in vorticity distribution inside the core.
Indeed, the second term on the right-hand side of the equation for \( Z_{i}\) is 
written in a more general form given by
\begin{align*}
\frac{\Gamma_{i}}{4\pi R_{i}}\left(
\log\frac{8R_{i}}{a_{i}}-\frac{1}{2}+\Delta(a_{i})\right), 
\end{align*}
where
\begin{align*}
\Delta(a_{i})=\frac{1}{\Gamma_{i}^{2}}\int^{a_{i}}_{0}
\frac{\gamma_{i}(s)^{2}}{s}{\rm d}s, \quad 
\gamma_{i}(s)=2\pi \int^{s}_{0}\omega_{i}(r)r {\rm d}r.
\end{align*}
Here, \( \gamma_{i}\) is the velocity circulation around the central part of the core and 
\( \omega_{i} \) is the vorticity distribution of the cross-section of the core.
The model considered in \cite{39} corresponds to vortex rings with a circular 
core cross-section with constant radius and a uniform vorticity distribution.
As was stated in the introduction, the ODE model has a long history behind its 
derivation and analysis, and is widely accepted as the model which describes the 
motion of interacting coaxial vortex rings.

On the other hand, we derived the model system (\ref{model}) as a system of
partial differential equations (PDE). The leapfrogging solutions obtained in 
Sections 3 and 4 are solutions of this PDE model, and hence it is 
possible to consider the stability (or the lack there of) of these solutions 
under non-symmetric perturbations, i.e. a perturbation in which the 
shape of the filament is deviated from a circle. This is not possible 
in the framework of the ODE model. 

By comparing the systems (\ref{ODE}) and (\ref{RZode}), one can say that
(\ref{ODE}) is more focused on the 
precise description of the motion of the rings. In fact, system 
(\ref{RZode}) can be seen as a simplification of the ODE model.
This is especially apparent for the term describing the self-induced velocity 
in the equations for \( z_{i}\) and \( Z_{i}\). 
This is expected since our model neglects the effects of the finite core size.
To this end, we give a numerical plot of the phase portraits of the system 
(\ref{2dyn}) to observe the possible dynamics of the circular filament pair 
described by our model.
We do this to show that although our model is simpler, it is still able to capture 
the essential characteristics of the possible motions of vortex rings. 
The following plots were obtained by Mathematica as the level-sets of the 
Hamiltonian \( {\cal H}(\theta,W)\). The parameters were set at
\( d=1\), \( \alpha =0.1\), and plots are given for 
\( \beta = 1,2,\) and \(4\).
\vspace*{5mm}
\begin{figure}[h]
\begin{center}
\begin{overpic}[width=0.6\textwidth]{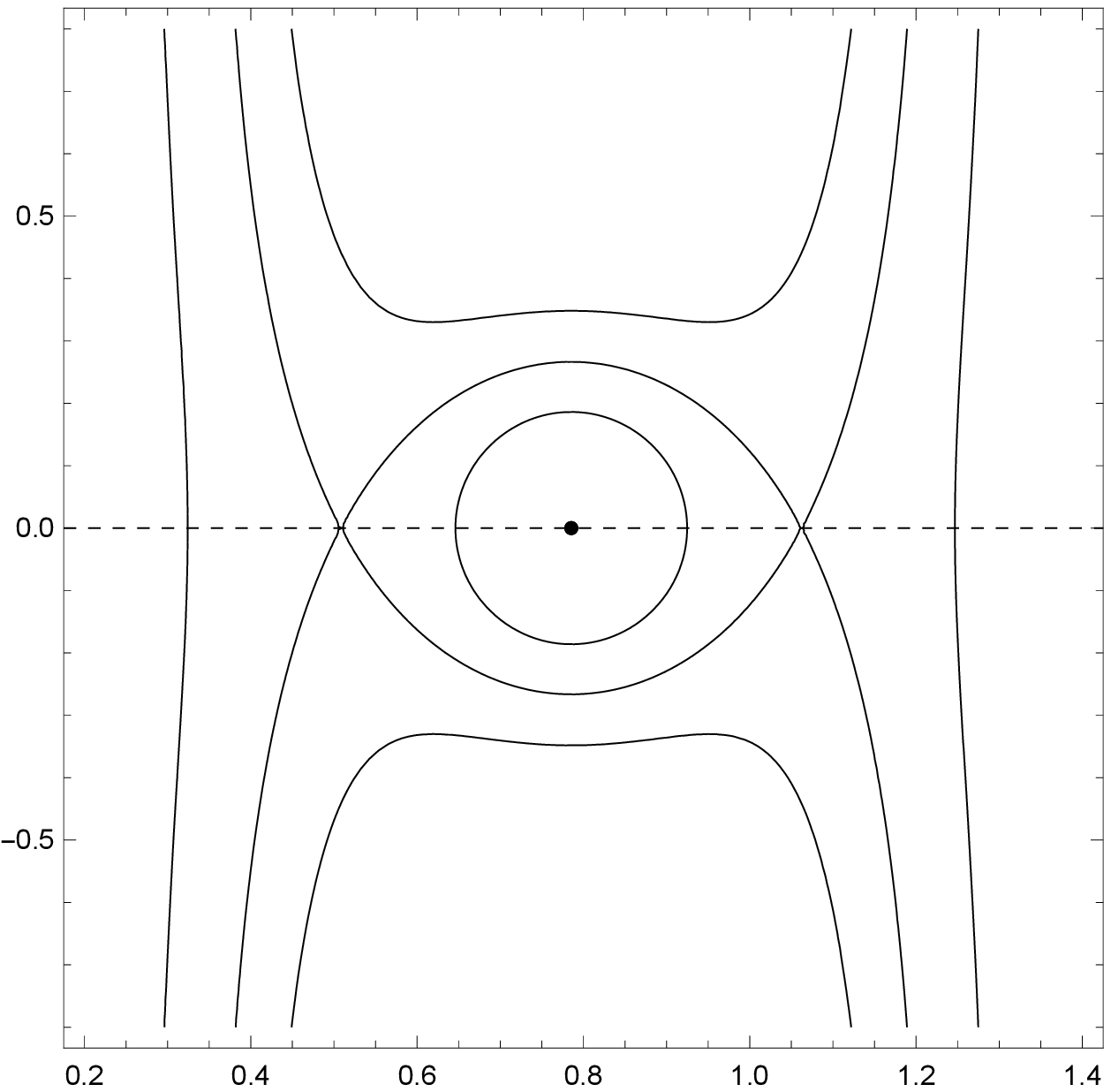}
\put(102,49){\( \theta\)}
\put(29,90){\( a\)}
\put(89,70){\( b\)}
\end{overpic}\\
\end{center}
\caption{\( \beta=1\)}
\end{figure}
\newpage
\begin{figure}[H]
\begin{center}
\begin{overpic}[width=0.6\textwidth]{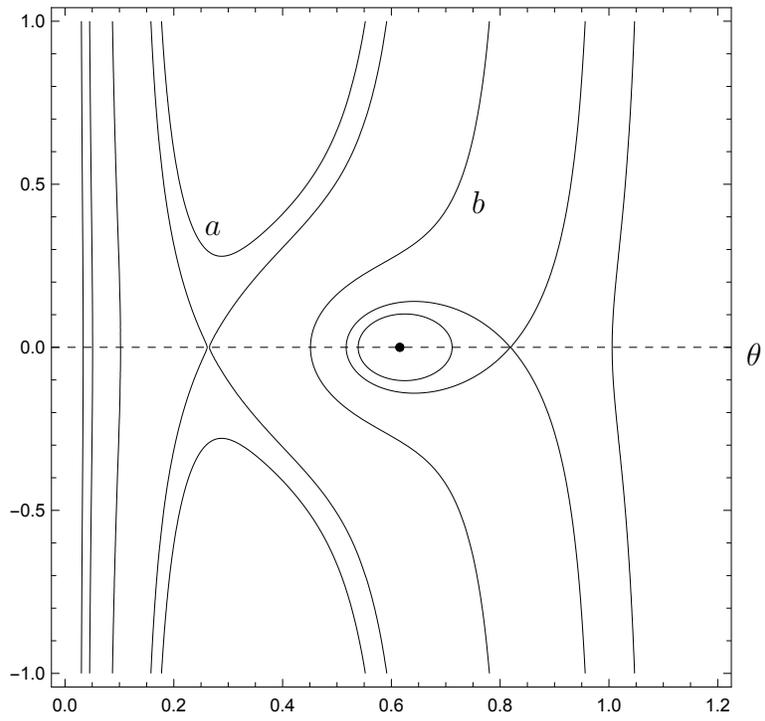}
\put(102,49){\( \theta\)}
\put(27,67){\( a\)}
\put(64,70){\( b\)}
\end{overpic} \\[5mm]
\begin{overpic}[width=0.6\textwidth]{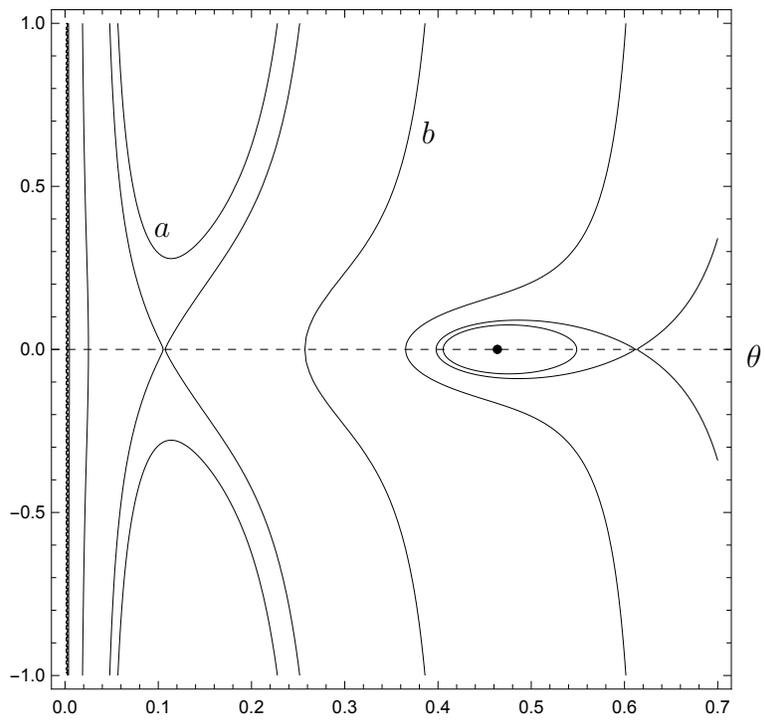}
\put(102,49){\( \theta\)}
\put(20,67){\( a\)}
\put(57,80){\( b\)}
\end{overpic}
\caption{top: \( \beta=2\), bottom: \( \beta=4\)}
\end{center}
\end{figure}
\newpage
\noindent
The dashed line is the \( \theta \)-axis, the vertical axis is the 
\( W \)-axis, and the black dot in each plot is the point \( (\theta_{\beta},0)\).
The phase portraits suggest that there are three possible types of motion:
leapfrogging, single passage, and repulsion. 

Leapfrogging motion are motions corresponding to closed orbits, which its existence was
rigorously proved in Section 3. We can see that in each of the three portraits, there are
closed orbits revolving around the point \( (\theta_{\beta},0)\).

Single passage is the motion in which one filament goes through the other 
once and then separate from each other indefinitely. The orbits labeled \( b\) correspond to
such motion.

Repulsion is the motion in which one filament approaches the other until some 
minimal distance is attained, and then is repulsed away from each other. Orbits 
labeled as \( a\) correspond to such motion. 

These dynamics qualitatively agree with the dynamics obtained in \cite{39}, 
supporting the validity of our model. 

\medskip

%
%
%
%
%
%
%
%
\subsection{Concluding Remarks}
We derived a PDE model describing the interaction of a pair of vortex filaments.
The system was explicitly solved for a pair of straight and parallel lines, 
and showed that the motion resembled that of a pair of point vortices.
We also proved rigorously that there exist solutions which correspond to 
leapfrogging motion of interacting coaxial circular filaments, and 
gave necessary and sufficient conditions for such motion to occur.

Although the model system (\ref{model}) is a system of partial differential equations,
the analysis carried out in this paper is essentially 
for systems of ordinary differential equations.
As the next step, the author would like to consider the unique solvability of 
(\ref{model}) in a mathematically rigorous setting, in a neighbourhood 
of exact solutions obtained in this paper at the very least, to 
investigate the stability of these solutions under non-symmetric perturbations.

\medskip

\section*{Acknowledgements}
This work was supported in part by JSPS Grant-in-Aid for 
Young Scientists(B) grant number 15K17579.

\vspace*{1cm}
\noindent
Masashi Aiki\\
Department of Mathematics\\
Faculty of Science and Technology, Tokyo University of Science\\
2641 Yamazaki, Noda, Chiba 278-8510, Japan\\
E-mail: aiki\verb|_|masashi\verb|@|ma.noda.tus.ac.jp

\end{document}